\documentclass[11pt]{amsart}
\usepackage{amsopn}
\usepackage{amsmath,amsthm,amssymb}
\newcommand{\nc}{\newcommand}

\nc{\n}{\noindent}
\nc{\Id}{\mathbf{1}}
\nc{\vs}{\vspace{8pt}}
\nc{\alt}{\raise1pt\hbox{$\bigwedge$}}
\nc{\ncp}{\nabla^\mathrm{CP}}
\nc{\nhc}{\nabla^\mathrm{HC}}
\nc{\wncp}{\widehat\nabla^\mathrm{CP}}
\nc{\ct}{\cos_\theta}
\nc{\st}{\sin_\theta}
\nc{\ctt}{\cos_{\theta/2}}
\nc{\stt}{\sin_{\theta/2}}
\nc{\ft}{\hbox{$\frac12$}}
\nc{\Aff}{\mathit{Aff}}

\nc{\vg}{\mathfrak{v} }
\nc{\wg}{\mathfrak{w} }
\nc{\zg}{\mathfrak{z} }
\nc{\ngo}{\mathfrak{n} }
\nc{\kg}{\mathfrak{k} }
\nc{\mg}{\mathfrak{m} }
\nc{\bg}{\mathfrak{b} }
\nc{\ggo}{\mathfrak{g} }
\nc{\ggob}{\overline{\mathfrak{g}} }
\nc{\sog}{\mathfrak{so} }
\nc{\sug}{\mathfrak{su} }
\nc{\spg}{\mathfrak{sp} }
\nc{\slg}{\mathfrak{sl} }
\nc{\glg}{\mathfrak{gl} }
\nc{\cg}{\mathfrak{c} }
\nc{\hg}{\mathfrak{h} }
\nc{\tg}{\mathfrak{t} }
\nc{\ug}{\mathfrak{u} }
\nc{\dg}{\mathfrak{d} }
\nc{\ag}{\mathfrak{a} }
\nc{\pg}{\mathfrak{p} }
\nc{\sg}{\mathfrak{s} }
\nc{\aff}{\mathfrak{aff}}
\nc{\pca}{\mathcal{P}}
\nc{\nca}{\mathcal{N}}
\nc{\vp}{\varphi}
\nc{\ddt}{\frac{{\rm d}}{{\rm d}t}}
\nc{\SO}{{\sf SO}}
\nc{\Spe}{{\sf Sp}}
\nc{\Sl}{{\sf Sl}}
\nc{\SU}{{\sf SU}}
\nc{\Or}{{\sf O}}
\nc{\U}{{\sf U}}
\nc{\Gl}{{\sf Gl}}
\nc{\Se}{{\sf S}}
\nc{\Cl}{{\sf Cl}}
\nc{\Spin}{{\sf Spin}}
\nc{\Pin}{{\sf Pin}}
\nc{\RR}{{\mathbb R}}
\nc{\HH}{{\mathbb H}}
\nc{\CC}{{\mathbb C}}
\nc{\ZZ}{{\mathbb Z}}
\nc{\FF}{{\mathbb F}}
\nc{\NN}{{\mathbb N}}
\nc{\GG}{{\mathbb G}}
\nc{\JJ}{{\mathbb J}}
\nc{\II}{{\mathbb I}}
\nc{\KK}{{\mathbb K}}
\nc{\DD}{{\mathbb D}}
\nc{\EE}{{\mathbb E}}
\nc{\VV}{{\mathbb V}}
\nc{\ad}{\operatorname{ad}}
\nc{\Ad}{\operatorname{Ad}}
\nc{\rank}{\operatorname{rank}}
\nc{\Irr}{\operatorname{Irr}}
\nc{\End}{\operatorname{End}}
\nc{\Aut}{\operatorname{Aut}}
\nc{\Inn}{\operatorname{Inn}}
\nc{\Der}{\operatorname{Der}}
\nc{\Ker}{\operatorname{Ker}}
\nc{\Iso}{\operatorname{I}}
\nc{\Le}{\operatorname{L}}
\nc{\tr}{\operatorname{tr}}
\nc{\dif}{\operatorname{d}}
\nc{\sen}{\operatorname{sen}}
\nc{\modu}{\operatorname{mod}}
\nc{\Ric}{\operatorname{R}}
\nc{\Sym}{\operatorname{Sym}}
\nc{\sca}{\operatorname{sc}}
\nc{\scalar}{{\sf s}}
\nc{\grad}{\operatorname{grad}}
\nc{\ricci}{\operatorname{r}}
\nc{\riccin}{\operatorname{Ric}}
\nc{\Lie}{\operatorname{L}}
\nc{\tang}{\operatorname{T}}
\nc{\Hom}{\operatorname{Hom}}
\theoremstyle{plain}
\newtheorem{thm}{Theorem}[section]
\newtheorem{prop}[thm]{Proposition}
\newtheorem{cor}[thm]{Corollary} 
\newtheorem{lem}[thm]{Lemma}
\theoremstyle{definition}
\newtheorem{defi}[thm]{Definition}
\theoremstyle{remark}
\newtheorem*{rem}{Remark}
\newtheorem*{rems}{Remarks}
\newtheorem{exam}[thm]{Example}
\newtheorem{exams}[thm]{Examples}

\newcommand{\ri}{{\rm (i)}}
\newcommand{\rii}{{\rm (ii)}}
\newcommand{\riii}{{\rm (iii)}}

\title[Complex product structures]
{\Large Complex product structures\\[3pt] on Lie algebras}

\author[A.~Andrada \,and \,S.~Salamon]{Adri\'an Andrada and Simon Salamon}

\begin{document}

\maketitle

\n\textbf{Abstract.} A study is made of real Lie algebras admitting compatible
complex and product structures, including numerous 4-dimensional examples. If
$\ggo$ is a Lie algebra with such a structure then its complexification has a
hypercomplex structure. It is shown in addition that $\ggo$ splits into the sum
of two left-symmetric subalgebras. Interpretations of these results are
obtained that are relevant to the theory of both hypercomplex and
hypersymplectic manifolds and their associated connections.

\smallbreak\n\textbf{MSC.} 17B60; 53C15, 53C30

\section{Introduction} 

A complex structure on a real Lie algebra $\ggo$ is an endomorphism $J$ of
$\ggo$ satisfying $J^2=-\Id$ and the usual integrability condition expressed in
terms of Lie brackets. Let $\ggo^{\CC}=\ggo\otimes_{\RR}\CC$ denote the
complexification of $\ggo$, and $\sigma:\ggo^\CC\to\ggo^\CC$ the corresponding
conjugation. The integrability of $J$ is then equivalent to a splitting
\[\ggo^{\CC}=\ggo^{1,0}\oplus\ggo^{0,1},\] where $\ggo^{1,0}$ and $\ggo^{0,1}$
are complex Lie subalgebras of $\ggo^\CC$ and
$\ggo^{0,1}=\sigma(\ggo^{1,0})$. A product structure on $\ggo$ is an
endomorphism $E$ of $\ggo$ which satisfies $E^2=\Id$ and an even simpler
integrability condition. Let \[\ggo=\ggo_+\oplus \ggo_-\] be the splitting in
which $\ggo_\pm$ denotes the eigenspaces corresponding to the eigenvalue $\pm
1$ of $E$. Then the integrability of $E$ is equivalent to requiring that
$\ggo_+,\,\ggo_-$ are both Lie subalgebras of $\ggo$.

In this article we are interested in another structure on $\ggo$ which arises
from the combination of the two above. Namely, a \textit{complex product
structure} on a real Lie algebra $\ggo$ is a pair $\{J,\,E\}$ where $J$ is a
complex structure on $\ggo$, $E$ is a product structure on $\ggo$, and
$JE=-EJ$. In this case, the subalgebras $\ggo_\pm$ satisfy $\ggo_-=J\ggo_+$.
We summarize the basic definitions in \S2, and provide some simple examples and
constructions of such structures. Whilst the concept of a product structure is
a relatively elementary one, we shall see that complex product structures are
altogether more subtle and their classification far from easy.

A complex product structure on a Lie algebra is an analogue of a hypercomplex
structure, i.e.\ a pair of anticommuting complex structures. Indeed, a complex
product structure is a `latent hypercomplex structure' in the sense that a
complex product structure on a Lie algebra $\ggo$ determines a hypercomplex
structure on the real Lie algebra $(\ggo^\CC)_\RR$ underlying $\ggo^\CC$. This
result, that we establish in \S3, provided our original interest in the
subject. We also explain how a standard complex product structure on
$\glg(2n,\RR)$ gives rise to a hypercomplex structure on $\glg(2n,\CC)$. In
general, $(\ggo^\CC)_\RR$ also admits a complex product structure, so the
process can be continued so as to obtain an infinite family of hypercomplex
structures.

Complex product structures are closely related to a number of other algebraic
concepts, and these relationships are pursued in \S4. If $(\ggo,J,E)$ is a Lie
algebra with a complex product structure as above, then $(\ggo,\ggo_+,\ggo_-)$
is an example of a double Lie algebra \cite{LW}. Actually, following
\cite{Maj}, $(\ggo_+,\,\ggo_-)$ is a matched pair of Lie algebras for certain
representations $\rho:\ggo_+\longrightarrow\glg(\ggo_-)$ and
$\mu:\ggo_-\longrightarrow\glg(\ggo_+)$. The isomorphism
$J:\ggo_+\longrightarrow\ggo_-$ then allows us to construct left-symmetric
algebra (LSA) structures on $\ggo_+$ and $\ggo_-$. This implies that neither
$\ggo_+$ nor $\ggo_-$ are semisimple \cite{B}, though we do not know whether
there are examples with $\ggo$ semisimple. Theorems \ref{lsa} and
\ref{converse} provide a complete characterization of Lie algebras carrying
complex product structures in terms of pairs of LSA structures, and we also
examine the special situation in which $\ggo$ is itself an LSA.

Whilst the focus of this paper is algebraic, the correspondence between flat
torsion-free connections on a Lie group and LSA structures on its Lie algebra
allows us to give various geometrical interpretations of the theory developed
thus far. This we do in \S5. Suppose that $\ggo$ has a complex product
structure, and is the Lie algebra of a Lie group $G$. The endomorphisms $J,E$
can be regarded as left-invariant tensors on $G$ still satisfying $JE=-EJ$.
This makes $G$ a \textit{complex product manifold}, a concept that is most
easily defined in terms of a $GL(n,\RR)$-structure for a diagonal inclusion
$GL(n,\RR)\subset GL(2n,\RR)$. We explain how properties of this type of
manifold reflect the theory of \S4, referring the reader to \cite{A} for a more
complete treatment.

A complex product manifold $M$ possesses two involutive totally real
distributions $T_\pm$. It also has a unique torsion-free connection $\ncp$ for
which $J,E$ are parallel, and the integral submanifolds of $T_\pm$ are flat
with respect to $\ncp$. Thus, complex product manifolds are integrable in a
rather strong sense. When $M=G$ is a Lie group, we obtain a local double Lie
group $(G,\,G_+,\,G_-)$ in which each of the subgroups $G_\pm$ has a flat
affine structure. We also show that $\ncp$ gives rise to the Obata connection
of the associated hypercomplex structure on $G^\CC$; indeed the two
connections share many properties.

There are natural links between complex product structures and symplectic
geometry, since a class of left-symmetric algebras are those arising from Lie
groups with a left-invariant symplectic form \cite{Boy}. It therefore makes
sense to look for complex product structures on Lie algebras for which
$\ggo_+,\ggo_-$ are both {\em symplectic} Lie algebras.  Algebraically, this
means that $n=2m$ is even and the structure group $GL(n,\RR)$ reduces to
$Sp(m,\RR)$. Geometrically, this arises when an associated Lie group is
endowed with a {\em hypersymplectic} structure in the sense of \cite{Hi}. We
explain this briefly in \S5, but do not pursue this aspect of the theory in the present
paper.

In \S6, we consider examples of 4-dimensional Lie algebras carrying complex
product structures, referring in part to the classification of complex structures given by Ovando 
\cite{O} and Snow \cite{Sn}. We determine the associated affine Lie groups and the hypercomplex
structures on their complexifications. We should emphasize though that the complex product
structures considered in this article can exist in dimensions $2n\not\equiv0$ mod~4, and we expect
a subsequent analysis of the case $n=3$ to prove fruitful.

\smallbreak\noindent{\small\textbf{Acknowledgments.} The first author was supported by
grants from CONICET, FONCYT and SECYT-UNC (Argentina). The authors are
grateful to I.~Dotti, A.~Fino and R.~Miatello for various valuable suggestions, and to
N.~Andruskiewitsch and D.~Burde for help with the references.}

\vs


\section{Preliminaries and basic results}

We begin by recalling some definitions which will be used throughout this
work. All Lie algebras will be finite dimensional and defined over $\RR$,
unless explicitly stated.

An {\em almost complex structure} on a Lie algebra $\ggo$ is a linear
endomorphism $J:\ggo \longrightarrow \ggo$ satisfying $J^2=-\Id$. If $J$
satisfies the condition
\begin{equation} \label{integrable}
J[X,Y]=[JX,Y]+[X,JY]+J[JX,JY] \quad \text{for all } X,Y\in \ggo,
\end{equation}
we will say that $J$ is {\em integrable} and we will call it a {\em complex
structure} on $\ggo$.  Note that the dimension of a Lie algebra carrying an
almost complex structure must be even. Given an almost complex structure $J$ on
$\ggo$, there is a splitting of $\ggo^{\CC}=\ggo\otimes_{\RR}\CC$ into the sum
of two subspaces
\begin{equation}\label{split}
 \ggo^{\CC}=\ggo^{1,0}\oplus\ggo^{0,1},
\end{equation}
where $\ggo^{1,0}$ and $\ggo^{0,1}$ are the eigenspaces corresponding to the
eigenvalues $\pm i$ of the complexification of $J$. It is easy to see that $J$
is integrable if and only if both $\ggo^{1,0}$ and $\ggo^{0,1}$ are complex
subalgebras of $\ggo^{\CC}$.

The splitting (\ref{split}) induces a decomposition of $(\ggo^\CC)^*$ as
\[ (\ggo^\CC)^*=(\ggo^{1,0})^*\oplus(\ggo^{0,1})^* \]
and then \[\alt^k(\ggo^\CC)^*=\bigoplus_{p+q=k}\Lambda^{p,q},\] where
$\Lambda^{p,q}\cong\alt^p(\ggo^{1,0})^*\otimes\alt^q(\ggo^{0,1})^*$ is the
space of $(p,q)$-forms relative to $J$.

Recall that the Lie bracket in $\ggo$ can be thought of as a linear map
$[\,,]:\alt^2 \ggo \longrightarrow \ggo$. In this way, we can consider its
transpose $\dif:\ggo^*\longrightarrow \alt^2\ggo^*$, defined as follows:
\[ (\dif f)(X\wedge Y)= -f([X,Y]) \quad \text{for all } f\in \ggo^*,\, 
X,Y\in \ggo.\] It is well known that $J$ is integrable if and only if
\[\dif(\Lambda^{1,0})\subseteq \Lambda^{2,0}\oplus\Lambda^{1,1}.\]

Whilst an almost complex structure on a Lie algebra $\ggo$ makes the latter a
complex vector space, $\ggo$ need not be a complex Lie algebra since the Lie
bracket might not be $\CC$-bilinear. This problem may be overcome in the
following situation. Let $\ggo$ be a Lie algebra with an almost complex
structure $I$ which satisfies
\begin{equation} \label{complexgroup} 
[IX,Y]=I[X,Y]\quad \text{for all }X,Y \in \ggo,
\end{equation}
It is straightforward to verify that $I$ is integrable, since
(\ref{complexgroup}) is stronger than (\ref{integrable}). Then $\ggo$ can be
viewed as a complex Lie algebra, with multiplication by $i$ given by the
endomorphism $I$, since (\ref{complexgroup}) implies the $\CC$-bilinearity of
the Lie bracket.

Condition (\ref{complexgroup}) is equivalent to $I\ad(X)=\ad(X)I$ for all $X\in
\ggo$, i.e.\ $\ad(X)$ is a complex transformation on $\ggo$, and it is also
equivalent to $\dif(\Lambda^{1,0})\subseteq \Lambda^{2,0}$. Conversely, if
$\hg$ is a complex Lie algebra, let $\hg_{\RR}$ denote the underlying real Lie
algebra. Then the endomorphism $I$ of $\hg_{\RR}$ given by multiplication by
$i$, i.e.\ $IX=iX$ for all $X\in\hg_{\RR}$, defines a complex structure on
$\hg_{\RR}$ satisfying (\ref{complexgroup}).

Next, we define another kind of structure on a Lie algebra which is analogous
to a complex structure. An {\em almost product structure} on $\ggo$ is a
linear endomorphism $E:\ggo \longrightarrow \ggo$ satisfying $E^2=\Id$ (and
not equal to $\pm\Id$). It is said to be {\em integrable} if \begin{equation}
\label{integrable2} E[X,Y]=[EX,Y]+[X,EY]-E[EX,EY] \quad \text{for all } X,Y\in
\ggo.  \end{equation} An integrable almost product structure will be called a
{\em product structure}. If $\dim \ggo_+=\dim\ggo_-$, where $\ggo_\pm$ is the
eigenspace of $\ggo$ associated to the eigenvalue $\pm1$ of $E$, then the
product structure $E$ is called a {\em paracomplex structure}
\cite{KK,L}. In this case, $\ggo$ also has even dimension.

Given an almost product structure $E$ on $\ggo$, we have a decomposition
$\ggo=\ggo_+\oplus\ggo_-$ into the direct sum of two linear subspaces, the
eigenspaces associated to $E$, which induces a splitting $\ggo^*=A_+\oplus A_-$
and then \[\alt^k\ggo^*=\bigoplus_{p+q=k}A^{p,q},\] where $A^{p,q}\cong\alt^p
A_+\otimes\alt^q A_-$. It is easy to verify that $E$ is integrable if and only
if $\ggo_+$ and $\ggo_-$ are both Lie subalgebras of $\ggo$, and this is in
turn equivalent to
\begin{equation} \label{decom}
\dif(A^{1,0})\subseteq A^{2,0}\oplus A^{1,1},\qquad\dif(A^{0,1})\subseteq
A^{0,2}\oplus A^{1,1}.
\end{equation}
Verification of these statements can be found within the proof of Proposition
\ref{equiv}.

Lie algebras carrying either a complex structure or a product structure have
many interesting properties that have been exhaustively studied. Our interest
lies in an appropriate combination of these two structures:

\vs

\begin{defi} \label{def1}
A {\em complex product structure} on the Lie algebra $\ggo$ is a pair
$\{J,\,E\}$ of a complex structure $J$ and a product structure $E$ satisfying
$JE=-EJ$.
\end{defi}

\vs

\n The condition $JE=-EJ$ implies that the eigenspaces corresponding to the
eigenvalues $+1$ and $-1$ of $E$ have the same dimension, showing that $E$ is
in fact a paracomplex structure on $\ggo$.

The endomorphism $F:=JE$ satisfies $F^2=\Id$, and overall $\{J,E,F\}$ obey the
rules 
\begin{equation}\label{rules}
\begin{array}{c} -J^2=E^2=F^2=\Id,\\[3pt]
JE=F,\quad EF=-J,\quad FJ=E,
\end{array} 
\end{equation} satisfied by the
so-called paraquaternionic numbers \cite{GMV}. It is easy to verify that
(\ref{integrable2}) is satisfied by $F$ in place of $E$. Indeed $\ggo$ has a
circle's worth $\{\ct E+\st F\}$ of product structures, and a corresponding
`pencil' of subalgebras $\ggo_\theta$ with $\ggo_0=\ggo_+$, 
\[ \ggo_{\pm\pi/2}=\{X\pm JX:X\in\ggo_+\},\] 
and $\ggo_\pi=\ggo_-$.

\vs

\begin{prop} \label{equiv}
Let $\ggo$ be a Lie algebra. The following statements are equivalent:
\begin{enumerate}
\item[\ri] $\ggo$ has a complex product structure,
\item[\rii] $\ggo$ has a complex structure $J$ and can be decomposed as $\ggo =
\ggo_+ \oplus \ggo_-$, where $\ggo_+,\ggo_-$ are Lie subalgebras of $\ggo$ and
$\ggo_-=J\ggo_+$.
\item[\riii] $\ggo$ has a complex structure $J$ and $\ggo^*$ can be decomposed
as $\ggo^* = A_+ \oplus A_-$, where $A_+,A_-$ are subspaces satisfying 
(\ref{decom}) and $A_-=JA_+$.
\end{enumerate}
\end{prop}

\begin{proof}
$\ri\Longleftrightarrow\rii$. If $\ggo$ has a complex product structure
$\{J,E\}$, let $\ggo_\pm$ denote the eigenspace corresponding to the eigenvalue
$\pm1$ of $E$. The integrability of $E$ implies that both $\ggo_+$ and $\ggo_-$
are Lie subalgebras of $\ggo$ and $JE=-EJ$ implies $\ggo_-=J\ggo_+$.

Conversely, if $\rii$ holds, set $E|_{\ggo_+}=\Id$ and $E|_{\ggo_-}=-\Id$. Then it
is easy to verify that $\{J,E\}$ determine a complex product structure on
$\ggo$.

\smallbreak\n $\rii\Longleftrightarrow\riii$. Assuming $\rii$, let
$A_+\subset\ggo^*$ (respectively, $A_-\subset\ggo^*$) denote the annihilator of
$\ggo_-$ (respectively, $\ggo_+$). Since $\ggo_+$ and $\ggo_-$ are subalgebras
of $\ggo$, one immediately checks (\ref{decom}). The complex structure $J$
induces a complex structure $J$ on the vector space $\ggo^*$ by setting
$Jf=f\circ J$. Then $\ggo_-=J\ggo_+$ implies $A_-=JA_+$.

Conversely, if $\riii$ holds, set $\ggo_+=\{X\in\ggo:f(X)=0 \text{ for all
}f\in A_-\}$ and $\ggo_-=\{X\in\ggo:f(X)=0 \text{ for all }f\in
A_+\}$. Conditions (\ref{decom}) imply at once that $\ggo_+$ and $\ggo_-$ are
Lie subalgebras of $\ggo$ and $\ggo_-=J\ggo_+$ follows from $A_-=JA_+$.
\end{proof}

\smallbreak

\begin{defi}
Let $(\ggo,J,E)$ and $(\ggo',J',E')$ be two Lie algebras with complex product
structures. We shall say that these complex product structures are {\em
equivalent} if there exists an isomorphism of Lie algebras $\phi:
\ggo\longrightarrow \ggo'$ such that $\phi J = J'\phi$ and $\phi E = E'\phi$.
\end{defi}

The equivalence of Lie algebras with complex product structures can be
reformulated in the following manner: according to $\rii$ of Proposition
\ref{equiv}, we have a decomposition into a sum of subalgebras
$\ggo=\ggo_+\oplus\ggo_-$ and $\ggo'=\ggo'_+\oplus\ggo'_-$ with
$\ggo_-=J\ggo_+,\,\ggo'_-=J'\ggo'_+$. Then $\phi E = E'\phi$ means that
$\phi(\ggo_{\pm})=\ggo'_{\pm}$, i.e.\ $\phi$ maps each eigenspace for $E$ onto
the corresponding eigenspace for $E'$. Also $\phi J = J'\phi$ means that the
natural extension of $\phi$ to an isomorphism
$\phi^{\CC}:\ggo^{\CC}\longrightarrow(\ggo')^{\CC}$ satisfies
$\phi^{\CC}(\ggo^{1,0})=(\ggo')^{1,0}$ and
$\phi^{\CC}(\ggo^{0,1})=(\ggo')^{0,1}$, where the splittings $
\ggo^{\CC}=\ggo^{1,0}\oplus\ggo^{0,1},\,
(\ggo')^{\CC}=(\ggo')^{1,0}\oplus(\ggo')^{0,1}$ are as in (\ref{split}).

A definition that will be useful for our purposes is the following, which
appeared in \cite{LW}.

\vs

\begin{defi} 
Three Lie algebras $(\ggo,\ggo_+,\ggo_-)$ form a {\em double Lie algebra} if
$\ggo_+$ and $\ggo_-$ are Lie subalgebras of $\ggo$ and
$\ggo=\ggo_+\oplus\ggo_-$ as vector spaces.
\end{defi}

\vs

Observe that a double Lie algebra $(\ggo,\ggo_+,\ggo_-)$ gives a product
structure $E:\ggo\longrightarrow\ggo$ on $\ggo$, where $E|_{\ggo_+}=\Id$ and
$E|_{\ggo_-}=-\Id$. Conversely, a product structure on the Lie algebra $\ggo$
gives rise to a double Lie algebra $(\ggo,\ggo_+,\ggo_-)$, where $\ggo_{\pm}$
is the eigenspace associated to the eigenvalue $\pm 1$ of $E$. Thus, a complex
product structure $\{J,E\}$ on $\ggo$ gives rise to a double Lie algebra
$(\ggo,\ggo_+\ggo_-)$ with $\ggo_-=J\ggo_+$. The complex structure $J$
determines a vector space isomorphism $J:\ggo_+\longrightarrow\ggo_-$ which
satisfies the condition (\ref{integrable}) on $\ggo$. In particular, $\dim
\ggo_+=\dim\ggo_-$.

Let us assume now we begin with a double Lie algebra $(\ggo,\ggo_+,\ggo_-)$
such that $\dim\ggo_+=\dim\ggo_-$. Consider a linear isomorphism
$\varphi:\ggo_+\longrightarrow\ggo_-$ and define an endomorphism $J$ of $\ggo$
by $J(X+A)=-\varphi^{-1}(A)+\varphi(X)$ for $X\in\ggo_+$ and $A\in\ggo_-$. $J$
is clearly an almost complex structure on $\ggo$ but, if we want $J$ to be
integrable, we need $\varphi$ to satisfy the extra condition
\begin{equation}\label{fi}
\varphi[X,Y]+\varphi^{-1}[\varphi X,\varphi Y]=[\varphi X,Y]+[X,\varphi Y],
\quad X,Y\in\ggo_+.
\end{equation}

These observations are summarized by

\vs

\begin{prop}\label{double}
A double Lie algebra $(\ggo,\ggo_+,\ggo_-)$ with $\dim\ggo_+=\dim\ggo_-$ is
associated to a complex product structure on the Lie algebra $\ggo$ if and only
if there exists a linear isomorphism $\varphi:\ggo_+\longrightarrow\ggo_-$ such
that (\ref{fi}) holds.
\end{prop}

\vs

We now wish to exhibit some simple examples of Lie algebras carrying complex
product structures. A natural class to begin with is the class of real Lie
algebras underlying a complex Lie algebra, since we have already mentioned that
they are naturally endowed with a complex structure. Nevertheless, it is
impossible to find interesting examples of complex product structures among
these algebras, due to

\vs

\begin{prop} \label{corcom}
Suppose that $\ggo$ admits a complex product structure given by a complex
structure $I$ and a product structure $E$ where $I$ satisfies
(\ref{complexgroup}). Then $\ggo$ is abelian.
\end{prop}

\begin{proof}
There is a splitting $\ggo = \ggo_+ \oplus \ggo_-$, with $\ggo_+$ and
$\ggo_-=I\ggo_+$ Lie subalgebras of $\ggo$. If $U,V \in\ggo_+$ then
\[ E[U,IV]=EI[U,V]=-IE[U,V]=-I[U,V]=-[U,IV],\]
which shows $[U,IV]\in\ggo_-$, and
\[ E[U,IV]=-E[I(IU),IV]=-EI[IU,IV]=IE[IU,IV]=-I[IU,IV]=[U,IV],\]
which shows $[U,IV]\in\ggo_+$. Thus, $\ggo_{\pm}$ are ideals in $\ggo$ and the
proposition follows from the following lemma.
\end{proof}

\smallbreak

\begin{lem}\label{ideales}
Let $\ggo$ admit a complex product structure $\{J,E\}$ and let
$(\ggo,\ggo_+,\ggo_-)$ be its associated Lie algebra. If both $\ggo_+$ and
$\ggo_-$ are ideals in $\ggo$ then $\ggo$ is abelian.
\end{lem}

\begin{proof}
There is a splitting $\ggo = \ggo_+ \oplus \ggo_-$, with $\ggo_+$ and
$\ggo_-=J\ggo_+$ ideals in $\ggo$. Since $[\ggo_+,\ggo_-]=0$, we have that
$[JU,V]=[U,JV]=0$ for all $U,V\in \ggo_+$. From the integrability of $J$ we
obtain $J[U,V]=J[JU,JV]$, and then $[U,V]=[JU,JV]=0$ for all $U,V\in
\ggo_+$. Hence $\ggo_{\pm}$ are abelian ideals and it follows that $\ggo$ is
abelian.
\end{proof}

\vs

\begin{exams}\label{ejemplo}
(i) Given any Lie algebra $\ug$, set $\ggo=\ug\times\ug$. Let
$E:\ggo\longrightarrow\ggo$ be defined by $E(U,V)=(U,-V)$. Then $E$ is clearly
a product structure on $\ggo$. However, Lemma \ref{ideales} implies that there
exists no complex product anticommuting with $E$ unless $\ug$ is abelian. In
this case any complex structure anticommuting with $E$ can be seen to be
equivalent to the standard complex structure on the abelian Lie algebra
$\CC^n=\RR^n\times\RR^n$ given by $J(U,V)=(-V,U)$.

\smallbreak\n(ii) The simplest non-abelian example of a Lie algebra carrying a
complex product structure is given by $\aff(\RR)$, the Lie algebra of the group
$\Aff(\RR)$ of affine motions of the real line. $\aff(\RR)$ is a solvable non
nilpotent Lie algebra with a basis $\{X,Y\}$ satisfying $[X,Y]=Y$. There is a
complex structure $J$ and a product structure $E$ on $\aff(\RR)$ given by
$JX=Y$ and $EX=X,\,EY=-Y$. As $JE=-EJ$, $\{J,E\}$ defines a complex product
structure on $\aff(\RR)$.
\end{exams}

\vs

An important class of Lie algebras carrying a complex product structure is
obtained by the following construction. Consider a finite dimensional real
vector space ${\mathcal A}$ equipped with a bilinear product ${\mathcal
A}\times{\mathcal A}\longrightarrow{\mathcal A},\; (a,b)\mapsto ab$. On the
vector space ${\mathcal A}\oplus{\mathcal A}$ we consider the following
skew-symmetric bilinear bracket:
\begin{equation}\label{bracket}
[(a,b),(a',b')]=(aa'-a'a,ab'-a'b), \qquad a,b,a',b'\in {\mathcal A}. 
\end{equation}
Suppose for a moment that the bracket (\ref{bracket}) satisfies the Jacobi
identity, hence defining a Lie algebra structure on ${\mathcal
A}\oplus{\mathcal A}$. Let $J:{\mathcal A}\oplus{\mathcal A}\longrightarrow
{\mathcal A}\oplus{\mathcal A}$ be the endomorphism given by
\[ J(a,b)=(-b,a), \qquad a,b\in {\mathcal A}.\]
A computation shows that $J$ is a complex structure on ${\mathcal
A}\oplus{\mathcal A}$. Furthermore, if we denote ${\mathcal A}_+:={\mathcal
A}\oplus\{0\}$ and ${\mathcal A}_-:=\{0\}\oplus {\mathcal A}$, then it is easy
to see that ${\mathcal A}_+$ is a subalgebra, ${\mathcal A}_-$ is an abelian
ideal of ${\mathcal A}\oplus{\mathcal A}$ and also ${\mathcal A}_-=J{\mathcal
A}_+$. Hence, letting $E:{\mathcal A}\oplus{\mathcal A}\longrightarrow
{\mathcal A}\oplus{\mathcal A}$ be given by $E|_{{\mathcal
A}_+}=\Id,\,E|_{{\mathcal A}_-}=-\Id$, we obtain a complex product structure
$\{J,E\}$ on ${\mathcal A}\oplus{\mathcal A}$.

Let us now investigate the conditions under which (\ref{bracket}) does in fact
satisfy the Jacobi identity. To do so, we shall need the following concept,
which has already been intensively studied (see for instance \cite{B1,B,H,Se}).

\vs

\begin{defi}
A {\em left-symmetric algebra (LSA)} structure on a Lie algebra $\hg$ is a
bilinear product $\hg\times\hg\longrightarrow\hg,\,(x,y)\mapsto x\cdot y$,
which satisfies
\begin{equation}\label{flat}
x\cdot(y\cdot z)-(x\cdot y)\cdot z=y\cdot(x\cdot z)-(y\cdot x)\cdot z
\end{equation}
and 
\begin{equation}\label{tfree} 
[x,y]=x\cdot y-y\cdot x.
\end{equation}
\end{defi}

\vs

\begin{lem}\label{JacLSA}
The bracket on ${\mathcal A}\oplus{\mathcal A}$ given by (\ref{bracket})
satisfies the Jacobi identity if and only if ${\mathcal A}$ is a Lie algebra
with an LSA structure.
\end{lem}

\begin{proof}
Let $a,a',a'',b,b',b''\in{\mathcal A}$. Then it is easy to see that
\[ [[(a,b),(a',b')],(a'',b'')]+[[(a',b'),(a'',b'')],(a,b)]+
[[(a'',b''),(a,b)],(a',b')]=0\] if and only if
\begin{multline*}
(aa')a''-a(a'a'')-(a'a)a''+a'(aa'')+(a''a)a'-a''(aa')-(aa'')a'+a(a''a')+ \\
(a'a'')a-a'(a''a)-(a''a')a+a''(a'a)=0
\end{multline*}
and 
\begin{multline}\label{abc}
(aa')b''-a(a'b'')-(a'a)b''+a'(ab'')+(a''a)b'-a''(ab')-(aa'')b'+a(a''b')+ \\
(a'a'')b-a'(a''b)-(a''a')b+a''(a'b)=0.
\end{multline}
It is clear that if ${\mathcal A}$ is a Lie algebra with an LSA structure, then
the Jacobi identity is valid on ${\mathcal A}\oplus{\mathcal A}$, using
(\ref{flat}). Conversely, suppose that the bracket (\ref{bracket}) satisfies
the Jacobi identity. If we take $a'=b''=0$ and $a=x,\,a''=y,\,b'=z$ in
(\ref{abc}), we obtain equation (\ref{flat}). As we have already mentioned,
${\mathcal A}_+={\mathcal A}\oplus\{0\}$ is a Lie subalgebra of ${\mathcal
A}\oplus{\mathcal A}$.  Identifying ${\mathcal A}_+$ with ${\mathcal A}$ in the
obvious way, ${\mathcal A}$ itself acquires a Lie algebra structure which
satisfies $[a,a']=aa'-a'a$. Hence ${\mathcal A}$ is a Lie algebra with an LSA
structure and the proof of the lemma is complete.
\end{proof}

\vs

\begin{rems} 
(i) LSA structures on a Lie algebra are also known as Koszul-Vinberg structures
or affine structures.

\smallbreak\n(ii) If we define $(x,y,z)=x\cdot(y\cdot z)-(x\cdot y)\cdot z$ for
all $x,y,z\in\hg$, condition (\ref{flat}) becomes $(x,y,z)=(y,x,z)$, whence the
name ``left-symmetric".

\smallbreak\n(iii) Lemma~\ref{JacLSA} generalizes the result given in
\cite{BD}, where it is shown that the bracket (\ref{bracket}) gives a Lie
algebra structure on $\aff({\mathcal A}):={\mathcal A}\oplus{\mathcal A}$
whenever ${\mathcal A}$ is an associative algebra. We will also use the
notation $\aff({\mathcal A})$ for the Lie algebra ${\mathcal A}\oplus{\mathcal
A}$ when ${\mathcal A}$ is a Lie algebra with an LSA structure.
\end{rems}

\vs

To sum up, starting with a Lie algebra with an LSA structure, we have
constructed a Lie algebra with a complex product structure in which one of the
eigenspaces corresponding to the product structure is an (abelian) ideal. We
now reveal the converse.

\begin{prop}\label{masideal}
Let ($\ggo,J,E$) be a Lie algebra with a complex product structure, and let
$(\ggo,\ggo_+,\ggo_-)$ be its associated double Lie algebra,
$\ggo_-=J\ggo_+$. Suppose $\ggo_-$ is an {\em ideal} in $\ggo$. Then $\ggo_-$
is abelian and hence $\ggo$ is isomorphic to the semidirect product
$\ggo_+\ltimes_{\ad_{\ggo}}\ggo_-$. Also, $\ggo_+$ carries an LSA structure,
which is given by
\[ X\cdot Y=-J[X,JY],\quad X,Y\in\ggo_+ \]
\end{prop}

\begin{proof}
Let $X,Y\in\ggo_-$. Since $J$ is integrable, we have
\begin{equation}\label{ideal}
J[X,Y]=[X,JY]+[JX,Y]+J[JX,JY].
\end{equation} 
Now, the left-hand side of (\ref{ideal}) is in $\ggo_+$, whereas the right-hand
side is in $\ggo_-$. So, both sides of (\ref{ideal}) are zero, showing that
$\ggo_-$ is an {\em abelian} ideal. Then $\ggo$ is isomorphic to the semidirect
product $\ggo_+\ltimes_{\rho}\ggo_-$, where the representation
$\rho:\ggo_+\longrightarrow\glg(\ggo_-)$ is simply given by
$\rho(X)=\ad(X),\,X\in\ggo_+$. To see that $\ggo_+$ carries an LSA structure,
we compute for $X,Y,Z\in\ggo_+$:
\begin{align*}
X\cdot Y-Y\cdot X & = -J[X,JY]+J[Y,JX] \\
                  & = -J([X,JY]+[JX,Y]) \\
                  & = -J(J[X,Y]-J[JX,JY])\\
                  & = [X,Y]
\end{align*}
using the integrability of $J$ and the fact that $\ggo_-$ is abelian. Thus,
(\ref{tfree}) holds. Also,
\begin{align*}
X\cdot(Y\cdot Z)-Y\cdot(X\cdot Z) & = -X\cdot (J[Y,JZ])+Y\cdot(J[X,JZ])\\
                                  & = -J[X,[Y,JZ]]+J[Y,[X,JZ]] \\
                                  & = J[JZ,[X,Y]] \qquad\text{(using Jacobi)}\\
                                  & = -J[[X,Y],JZ] \\
                                  & = [X,Y]\cdot Z,       
\end{align*}
from where (\ref{flat}) follows. 
\end{proof}

\vs

\section{Hypercomplex structures from complex product structures}

In this section we will show that each complex product structure on a Lie
algebra gives rise to a hypercomplex structure on the real Lie algebra
underlying its complexification. Firstly, we recall

\vs

\begin{defi}
A {\em hypercomplex structure} on the Lie algebra $\ggo$ is a pair
$\{J_1,\,J_2\}$ of complex structures on $\ggo$ satisfying $J_1J_2=-J_2J_1$.
\end{defi}

\n Note that $J_3:=J_1J_2$ is another complex structure on $\ggo$ and
$\{J_1,J_2,J_3\}$ satisfy
\[ J_1J_2=J_3,\quad J_2J_3=J_1,\quad J_3J_1=J_2.\]\smallbreak 

Throughout this section, we shall use the notation: $\hat{\ggo}=
(\ggo^\CC)_\RR$. We have already noted that multiplication by $i$ in
$\ggo^\CC=\ggo\otimes_\RR\CC$ defines a complex structure $I$ on $\hat\ggo$
satisfying (\ref{complexgroup}). Now suppose that $\ggo$ has dimension $2n$ and
a complex structure $J$. In this situation, the complexification of $J$ can be
regarded as a complex structure $\JJ$ on $\hat\ggo$ that commutes with $I$. To
see this, use the decomposition
\[ \hat{\ggo}=\ggo \oplus I\ggo, \]
set $\JJ|_{\ggo}=J$ and extend the definition of $J$ to $I\ggo$ by
$\JJ(IX)=I(JX)$. It is easy to verify the integrability of $\JJ$, using
(\ref{complexgroup}).

Of course, $I$ and $\JJ$ never determine a hypercomplex structure because they
commute. Our aim is to construct hypercomplex structures on $\ggo$ by retaining
$\JJ$, but modifying $I$.

\vs

\begin{lem}\label{III}
Let $\ggo$ be a Lie algebra with a complex structure $I$ satisfying
(\ref{complexgroup}). Suppose there is a splitting $\ggo=\ug_1\oplus\ug_2$ with
$\ug_1,\,\ug_2$ complex subalgebras of $\ggo$.  Then the linear endomorphism
$\II$ defined by
\[\II|_{\ug_1}=I,\quad\II|_{\ug_2}=-I\]
is a complex structure on $\ggo$.
\end{lem}

\begin{proof}
It is clear that $\II^2=-\Id$, so, we only have to check the integrability of
$\II$. Let $X=X_1+X_2\in\ggo,\,Y=Y_1+Y_2\in\ggo$, with $X_1,Y_1\in\ug_1$ and
$X_2,Y_2\in\ug_2$. Then we have
\begin{align*}
\II[X,Y]& =\II[X_1+X_2,Y_1+Y_2] \\
        & =I[X_1,Y_1]-I[X_2,Y_2]+\II([X_1,Y_2]+[X_2,Y_1])
\end{align*}

On the other hand,
\begin{align*}
[\II X,Y]& =[IX_1-IX_2,Y_1+Y_2]\\
         & =I[X_1,Y_1]+I[X_1,Y_2]-I[X_2,Y_1]-I[X_2,Y_2]\\
[X,\II Y]& =[X_1+X_2,IY_1-IY_2]\\
         & =I[X_1,Y_1]-I[X_1,Y_2]+I[X_2,Y_1]-I[X_2,Y_2]\\
\II[\II X,\II Y]& =\II[IX_1-IX_2,IY_1-IY_2]\\
                & =\II(-[X_1,Y_1]+[X_1,Y_2]+[X_2,Y_1]-[X_2,Y_2]\\
                & =-I[X_1,Y_1]+I[X_2,Y_2]+\II([X_1,Y_2]+[X_2,Y_1])
\end{align*}
Hence, $\II[X,Y]=[\II X,Y]+[X,\II Y]+\II[\II X,\II Y]$ and $\II$ is integrable.
\end{proof}

\vs

\begin{thm} \label{hyper}
If $\ggo$ has a complex product structure then $\hat{\ggo}$ has a hypercomplex
structure $\{\II,\JJ\}$ with $\II$ and $\JJ$ as above.
\end{thm}

\begin{proof} From Proposition \ref{equiv}, we have a decomposition 
$\ggo=\ggo_+\oplus\ggo_-$, with $\ggo_+$ and $\ggo_-=J\ggo_+$ subalgebras of
$\ggo$. By complexifying, we obtain $\hat{\ggo}=(\ggo_+\oplus I\ggo_+)\oplus
(\ggo_-\oplus I\ggo_-)$, where $\ggo_+\oplus I\ggo_+$ and $\ggo_-\oplus
I\ggo_-$ are {\em complex} subalgebras of $\hat{\ggo}$. By the lemma, the
endomorphism $\II$ defined by $\II|_{\ggo_+\oplus I\ggo_+}=I,
\quad\II|_{\ggo_-\oplus I\ggo_-}=-I$ is a complex structure on
$\hat{\ggo}$. Also, it is clear that $\II$ and $\JJ$ anticommute.
\end{proof}

\vs

\begin{rem}
It is easy to verify that two equivalent complex product structures on the Lie
algebra $\ggo$ give rise to equivalent hypercomplex structures on
$\hat{\ggo}$. We need only recall that two hypercomplex structures
$\{J_k\}_{k=1,2}$ and $\{J'_k\}_{k=1,2}$ are said to be equivalent if there
exists an automorphism $\phi$ of $\ggo$ such that $\phi J_k=J'_k\phi$ for
$k=1,2$.
\end{rem}

\vs

Note that $\hat{\ggo}$ also has a complex product structure, induced by the one
on $\ggo$. Indeed, define $\EE:\hat{\ggo}\longrightarrow\hat{\ggo}$ by
$\EE=-I\II$. As $I$ commutes with $\II$, we have that $\EE^2=\Id$ and
$\JJ\EE=-\EE\JJ$. It is also easy to verify that $\EE$ is integrable, as a
consequence of (\ref{complexgroup}). Finally, the eigenspaces corresponding to
the eigenvalues $+1$ and $-1$ of $\EE$ are, respectively, the subalgebras
$\ggo_+\oplus I\ggo_+$ and $\ggo_-\oplus I\ggo_-$. Hence, the pair
$\{\JJ,\EE\}$ determines a complex product structure on $\hat{\ggo}$.  Applying
Theorem \ref{hyper} to $\hat{\ggo}$, we obtain another Lie algebra with a
hypercomplex structure, with dimension twice the dimension of $\hat{\ggo}$ and
four times the dimension of $\ggo$. In this way, starting with a
$2n$-dimensional Lie algebra $\ggo$ endowed with a complex product structure,
we obtain a family $\{\ggo_{(k)}:k\in\NN\}$ of Lie algebras carrying
hypercomplex structures, where $\ggo_{(1)}:=\ggo,\,
\ggo_{(k+1)}:=\widehat{\ggo_{(k)}}$ and then $\dim\ggo_{(k)}=2^{k+1}n$.

\vs

\begin{exams} (i) If $\ggo=\aff(\RR)$, applying Theorem \ref{hyper} we 
obtain a hypercomplex structure on $\hat{\ggo}=\aff(\CC)$, the Lie algebra of
the group $\Aff(\CC)$ of affine motions of the complex plane. $\aff(\CC)$ has a
basis $\{X,Y,\hat{X},\hat{Y}\}$ with non-zero bracket relations
\[ [X,Y]=-[\hat{X},\hat{Y}]=Y,\quad [X,\hat{Y}]=[\hat{X},Y]=\hat{Y},\] 
and the hypercomplex structure $\{\II,\JJ\}$ is given by
\[ \II X=\hat{X},\quad \II Y=-\hat{Y},\quad \II^2=-\Id, \]
\[ \JJ X=Y,\quad \JJ \hat{X}=\hat{Y},\quad \JJ^2=-\Id. \] 
According to \cite{Bar}, $\aff(\CC)$ is the only 4-dimensional Lie algebra with
2-dimensional commutator ideal admitting a hypercomplex structure.

\smallbreak\n(ii) More generally, recall from \S2 that $\aff({\mathcal A})$
admits a complex product structure, where ${\mathcal A}$ is a Lie algebra with
an LSA structure. Hence, Theorem \ref{hyper} gives a hypercomplex structure on
$\widehat{\aff({\mathcal A})}=\aff({\mathcal A}^{\CC})$, where ${\mathcal
A}^{\CC}$ is the complex Lie algebra with an LSA structure obtained by
complexifying ${\mathcal A}$. This hypercomplex structure $\{\II,\JJ\}$ is
given by
\[ \II(a,b)=(ia,-ib),\qquad \JJ(a,b)=(-b,a), \]
where $a,b\in{\mathcal A}^{\CC}$. This hypercomplex structure coincides (up to
a sign) with the hypercomplex structure constructed in \cite{BD} (where only
associative algebras were considered).
\end{exams}

\vs

We will prove next that the reductive Lie algebra $\glg(2n,\RR)$ carries a
canonical complex product structure; hence we will obtain a hypercomplex
structure on $\glg(2n,\CC)$.

\vs

\begin{prop}\label{gl}
$\glg(2n,\RR)$ has a complex product structure for all $n\ge1$.
\end{prop}

\begin{proof}
First let
\begin{equation}\label{E0J0} 
E_0=\begin{pmatrix}
\Id & 0 \cr 
0 & -\Id \cr
\end{pmatrix} \quad 
\text{ and } \quad J_0=\begin{pmatrix} 0 & -\Id \cr
\Id & 0 \cr \end{pmatrix}
\end{equation}
denote the standard almost product and complex structures acting on $\RR^{2n}$,
where $\Id$ is the $(n\times n)$-identity matrix. Here $E_0$ and
$J_0$ are expressed as $2n\times 2n$ matrices acting by left multiplication on
column vectors. Then $E_0$ and $J_0$ define an almost product structure $E$ and
an almost complex structure $J$ on the set $\glg(2n,\RR)$ of $2n\times 2n$
matrices by right multiplication:
\[ E(A)=AE_0,\qquad J(A)=AJ_0 \]
for any $A\in\glg(2n,\RR)$. Clearly, $EJ=-JE$. To see that $E$ is integrable,
one can argue as follows.

Let $V=\RR^{2n}$, so $\glg(2n,\RR)=\End(V)=V\otimes V^*$. Multiplying two
matrices corresponds to
$(v\otimes\alpha)\cdot(w\otimes\beta)=\alpha(w)(v\otimes\beta)$ for
$\alpha,\beta\in V^*$ and $v,w\in V$. On a generator $v\otimes\alpha\in
V\otimes V^*$, the action of $E$ is given by $E(v\otimes\alpha)=v\otimes
E_0\alpha$, where $E_0\alpha=\alpha\circ E_0$. So, if $V_\pm^*$ denotes the
$(\pm1)$-eigenspace of $E_0$ on $V^*$, then the $(+1)$- and $(-1)$-eigenspaces
of $E$ are given, respectively, by $V\otimes V^*_+$ and $V\otimes V^*_-$. It
is easily seen that the product (and thus the Lie bracket) of two elements of
$V\otimes V^*_+$ will produce a matrix that remains in that subspace, showing
that $V\otimes V^*_+$ is a Lie subalgebra of $\End V$. The same is true for
$V\otimes V^*_-$ and hence $E$ is integrable.

In the same way we may prove that $J$ is integrable, by considering the $(\pm
i)$-eigenspaces of $J^{\CC}_0$ on $V^\CC$. Therefore, $\{J,E\}$ is a complex
product structure on $\glg(2n,\RR)$.
\end{proof}

\begin{cor}
The induced hypercomplex structure on $\glg(2n,\CC)$ is given by right
multiplication by the matrices
\[\II=\begin{pmatrix}
               i\Id & 0 \cr
               0 & -i\Id \cr
             \end{pmatrix},\qquad\JJ=J_0.\]
\end{cor}

\n This hypercomplex structure on $\glg(2n,\CC)$ merely corresponds to writing
\begin{equation}\label{2n}\textstyle \glg(2n,\CC)=(\CC^{2n})^*\otimes\HH^n
\ \cong\ \bigoplus\limits_{2n} \HH^n,\end{equation} and using right
multiplication by unit quaternions on $\HH^n$.

\vs

We close this section with a negative result, namely the failure of a na\"\i ve
attempt to extend Proposition~\ref{gl} to the Lie algebra
\begin{equation}\label{sp} \spg(n,\RR)=\{X\in\glg(2n,\RR):X^tJ_0+J_0X=0\}
\end{equation} of the symplectic group preserving a non-degenerate skew form 
on $\RR^{2n}$ (notation of (\ref{E0J0})).

\vs

\begin{exam} It can be seen that $X\in\spg(n,\RR)$ if and only if \[
X=\begin{pmatrix}A & B \cr C & -A^t \cr \end{pmatrix},\] where
$A,B,C\in\glg(n,\RR)$ and $B$ and $C$ are symmetric. The subspace of
$\spg(n,\RR)$ given by 
\[ \left\{\begin{pmatrix} A & 0 \cr 0 & -A^t
\end{pmatrix}: A\in\glg(n,\RR)\right \} \] 
is in fact a subalgebra of
$\spg(n,\RR)$ isomorphic to $\glg(n,\RR)$. Also, the subspaces 
\[\ag_+=\left\{\begin{pmatrix} 0 & B \cr 0 & 0
\end{pmatrix}:\,B\text{ symmetric}\right\},\qquad
\ag_-=\left\{\begin{pmatrix} 0 & 0 \cr C & 0
\end{pmatrix}:\,C\text{ symmetric}\right \}\] 
are abelian subalgebras of $\spg(n,\RR)$. We have thus obtained a decomposition 
\[ \spg(n,\RR)=\glg(n,\RR)\oplus\ag_+\oplus\ag_-\] 
of (\ref{sp}) into the direct sum of three subalgebras. The proof of the following result is based
on simple matrix computations best left to the reader.

\vs

\begin{prop}
There is no complex product structure $\{J,E\}$ on $\spg(n,\RR)$ such that
$J\ag_+=\ag_-$.
\end{prop}
\end{exam}

\vs

\section{Complex product structures from matched pairs of Lie algebras}
\label{pairs}

In this section we will characterize Lie algebras carrying a complex product
structure in terms of double Lie algebras and matched pairs of Lie algebras
endowed with a left-symmetric algebra (LSA) structure.

In Proposition \ref{double}, we characterized those double Lie algebras
associated to a complex product structure. Given a double Lie algebra
$(\ggo,\ggo_+,\ggo_-)$, the existence of such a structure on $\ggo$ depends on
the existence of a linear isomorphism $\varphi:\ggo_+\longrightarrow\ggo_-$
which satisfies a certain property that involves the Lie bracket of $\ggo$, not
only the brackets of $\ggo_+$ and $\ggo_-$. On the other hand, we should like
to construct Lie algebras with a complex product structure beginning with two
Lie algebras of the same dimension, which are not a priori subalgebras of
another one. The idea is to obtain certain conditions on these two Lie algebras
which ensure that their direct sum as vector spaces admits a Lie bracket with
respect to which both summands are Lie subalgebras.

Firstly, let us recall the following construction, which appears for example in
\cite{Maj,Mas}. Let ($\ug,\vg,\rho,\mu$) be a {\em matched pair of Lie
algebras}, i.e.\ $\ug$ and $\vg$ are Lie algebras and
$\rho:\ug\longrightarrow\glg(\vg)$ and $\mu:\vg\longrightarrow\glg(\ug)$ are
representations satisfying
\begin{gather} 
\rho(X)[A,B]-[\rho(X)A,B]-[A,\rho(X)B]+\rho(\mu(A)X)B-\rho(\mu(B)X)A=0,
\label{jacobi1}\\
\mu(A)[X,Y]-[\mu(A)X,Y]-[X,\mu(A)Y]+\mu(\rho(X)A)Y-\mu(\rho(Y)A)X=0,\label{jacobi2}
\end{gather}
for $X,Y\in\ug$ and $A,B\in\vg$. In this case, the vector space
$\ggo=\ug\oplus\vg$ admits a Lie bracket, given by
\[ [(X,A),(Y,B)]=([X,Y]+\mu(A)Y-\mu(B)X,[A,B]+\rho(X)B-\rho(Y)A),\]
for $X,Y\in\ug$ and $A,B\in\vg$. We will denote this new Lie algebra by
$\ggo=\ug \bowtie^{\rho}_{\mu}\vg$ (or simply $\ggo=\ug \bowtie\vg$) and will
call it the {\em bicrossproduct} of $\ug$ and $\vg$. Observe that
$\ug\equiv\ug\oplus\{0\}$ and $\vg\equiv\{0\}\oplus\vg$ are subalgebras of
$\ggo$ and then $(\ug\bowtie\vg,\ug,\vg)$ is a double Lie algebra.

Conversely, if $\ug$ and $\vg$ are Lie subalgebras of a Lie algebra $\ggo$ such
that $\ug\oplus\vg=\ggo$, i.e.\ $(\ggo,\ug,\vg)$ is a double Lie algebra, then
($\ug,\vg,\rho,\mu$) forms a matched pair of Lie algebras, where the
representations $\rho:\ug\longrightarrow\glg(\vg)$ and
$\mu:\vg\longrightarrow\glg(\ug)$ are determined by
\begin{equation}\label{mu}
[X,A]=-\mu(A)X+\rho(X)A,\qquad X\in\ug,\,A\in\vg, 
\end{equation}
so that $\ggo=\ug\bowtie\vg$. 

\smallbreak

Let us now suppose that $\ggo$ is a Lie algebra endowed with a complex product
structure $\{J,E\}$ and let $(\ggo,\ggo_+,\ggo_-)$ be its associated double Lie
algebra, where $\ggo_-=J\ggo_+$. From the paragraph above, there exist
representations $\rho:\ggo_+\longrightarrow\glg(\ggo_-)$ and
$\mu:\ggo_-\longrightarrow\glg(\ggo_+)$ such that $(\ggo_+,\ggo_-,\rho,\mu)$
forms a matched pair of Lie algebras and then $\ggo=\ggo_+\bowtie\ggo_-$. In
this situation, a particular feature is given by the existence of the
isomorphism $J:\ggo_+\longrightarrow\ggo_-$, since it allows us to produce a
pair of representations equivalent to $\rho$ and $\mu$ in the following
way. Let $\tilde{\rho}:\ggo_+\longrightarrow\glg(\ggo_+)$ and
$\tilde{\mu}:\ggo_-\longrightarrow\glg(\ggo_-)$ be given by
\begin{equation}\label{tilde} 
\tilde{\rho}(X):=-J\rho(X)J,\qquad \tilde{\mu}(A):=-J\mu(A)J, 
\end{equation}
for $X\in\ggo_+$ and $A\in\ggo_-$. It is a simple matter to check that
$\tilde{\rho}$ and $\tilde{\mu}$ are indeed representations of $\ggo_+$ and
$\ggo_-$, respectively. With this notation, (\ref{mu}) translates into
\[ [X,A]=-J\tilde{\rho}(X)JA+J\tilde{\mu}(A)JX, \]
for $X\in\ggo_+$ and $A\in\ggo_-$. So, if $\pi_+:\ggo\longrightarrow\ggo_+$ and
$\pi_-:\ggo\longrightarrow\ggo_-$ denote the projections, we have the relations
\[ \tilde{\rho}(X)Y=-\pi_+J[X,JY],\qquad \tilde{\mu}(A)B=-\pi_-J[A,JB], \]
for $X,Y\in\ggo_+$ and $A,B\in\ggo_-$.

The relationship between Lie algebras with a complex product structure and Lie
algebras with an LSA is the content of the following result, which is a
generalization of Proposition \ref{masideal}.

\begin{thm}\label{lsa}
Let $(\ggo,J,E)$ be a Lie algebra with a complex product structure and let
$(\ggo,\ggo_+,\ggo_-)$ be its associated double Lie algebra with
$\ggo_-=J\ggo_+$. Then $\ggo_+$ and $\ggo_-$ carry an LSA structure.
\end{thm}

\begin{proof}
Define a bilinear product on $\ggo_+$ in the following manner
\[ \ggo_+\times\ggo_+\longrightarrow\ggo_+, \quad 
(X,Y)\mapsto X\cdot Y:=\tilde{\rho}(X)Y,\] where $\tilde{\rho}$ is as in
(\ref{tilde}). Let us first verify (\ref{tfree}). For $X,Y\in\ggo_+$, we have
\begin{align*}
X\cdot Y-Y\cdot X & = \tilde{\rho}(X)Y-\tilde{\rho}(Y)X \\
                  & = -\pi_+J([X,JY]+[JX,Y]) \\
                  & = -\pi_+J(J[X,Y]-J[JX,JY])\\
                  & = \pi_+([X,Y]-[JX,JY])\\
                  & = [X,Y]. 
\end{align*}
The third equality follows from the integrability of $J$. Next, to verify
(\ref{flat}), we perform the following computation. For $X,Y,Z\in\ggo_+$,
\[ X\cdot(Y\cdot Z)-Y\cdot(X\cdot Z)=\tilde{\rho}(X)\tilde{\rho}(Y)Z-
\tilde{\rho}(Y)\tilde{\rho}(X)Z=\tilde{\rho}([X,Y])Z=[X,Y]\cdot Z \]
since $\tilde{\rho}$ is a representation. Also,
\[ (X\cdot Y)\cdot Z-(Y\cdot X)\cdot Z=(X\cdot Y-Y\cdot X)Z=[X,Y]\cdot Z. \] 
Hence, $X\cdot(Y\cdot Z)-Y\cdot(X\cdot Z)=(X\cdot Y)\cdot Z-(Y\cdot X)\cdot Z$
and this implies (\ref{flat}).

To show that $\ggo_-$ carries an LSA structure too, consider the bilinear
product
\[ \ggo_-\times\ggo_-\longrightarrow\ggo_-, \quad 
(A,B)\mapsto A\cdot B:=\tilde{\mu}(A)B,\] where $\tilde{\mu}$ is as in
(\ref{tilde}). In the same way as before, it is shown that this product
satisfies (\ref{flat}) and (\ref{tfree}).
\end{proof}

\smallskip

\begin{rems} (i) A semisimple Lie algebra does not admit any LSA structure
(see \cite{B,H}). Thus $\ggo_+$ and $\ggo_-$ cannot be semisimple.

\smallbreak\n(ii) There do exist compact Lie algebras with an LSA
structure. The simplest is $\sog(3)\oplus\RR$, the Lie algebra of the Lie
groups $U(2)$ and $S^3\times S^1$. An LSA
structure on $\sog(3)\oplus\RR$ is obtained from the associative algebra $\HH$ of
quaternions. (This is related to the fact that $S^3\times S^1$ is an affine quotient of
$\HH\setminus\{0\}$, and admits a left-invariant flat torsion-free connection
\cite{Bar1}; see \S5).

\smallbreak\n(ii) Many, but not all, nilpotent Lie algebras admit LSA
structures \cite{Ben}.  \end{rems}

\vs

Conversely, we shall show that Lie algebras arising from certain pairs of Lie
algebras with LSA structures, admit a complex product structure.

\begin{thm} \label{converse}
Let $\ug$ and $\vg$ be two Lie algebras of the same dimension, both carrying an
LSA structure.  Suppose there exists a linear isomorphism
$\varphi:\ug\longrightarrow \vg$ such that the representations
$\rho:\ug\longrightarrow\glg(\vg)$ and $\mu:\vg\longrightarrow\glg(\ug)$ given
by
\[ \rho(X)A=\varphi(X\cdot \varphi^{-1}(A)),\qquad \mu(A)X=\varphi^{-1}(A\cdot\varphi(X)) \] 
satisfy (\ref{jacobi1}) and (\ref{jacobi2}). Then $(\ug,\vg,\rho,\mu)$ forms a
matched product of Lie algebras and the bicrossproduct Lie algebra
$\ggo=\ug\bowtie\vg$ admits a complex product structure.
\end{thm}

\begin{proof}
We start by defining an endomorphism $E:\ggo\longrightarrow\ggo$ by
$E|_{\ug}=\Id,\,E|_{\vg}=-\Id$.  Then $E^2=\Id$ and since the corresponding
eigenspaces are Lie subalgebras of $\ggo$, $E$ determines a product structure
on $\ggo$.  Next, we define another endomorphism $J$ of $\ggo$ by
\[ J(X,A)=(-\varphi^{-1}(A),\varphi(X)),\quad X\in\ug,\,A\in\vg.\]
It is clear that $J^2=-\Id$. We only have to prove the integrability. In order to
do so, we will consider several cases.

\n(i) For $X,Y\in\ug$, we have:
\[ J[(X,0),(Y,0)]=J([X,Y],0)=(0,\varphi[X,Y]). \]
On the other hand,
\begin{align*}
[J(X,0),(Y,0)] & = [(0,\varphi X), (Y,0)] \\
               & = -(-\mu(\varphi X)Y,\rho(Y)\varphi X)\\
               & = (\varphi^{-1}(\varphi X \cdot \varphi Y),-\varphi(Y\cdot X)),\\
[(X,0),J(Y,0)] & = [(X,0),(0,\varphi Y)] \\
               & = (-\mu(\varphi Y)X,\rho(X)\varphi Y)\\
               & = (-\varphi^{-1}(\varphi Y \cdot \varphi X),\varphi(X\cdot Y))
\end{align*}
and 
\begin{align*}
J[J(X,0),J(Y,0)] & = J[(0,\varphi X), (0,\varphi Y)] \\
                 & = J(0,[\varphi X,\varphi Y]) \\
                 & = (-\varphi^{-1}[\varphi X,\varphi Y],0).
\end{align*}
Therefore, using (\ref{tfree}), we obtain 
\[ J[(X,0),(Y,0)]=[J(X,0),(Y,0)]+[(X,0),J(Y,0)]+J[J(X,0),J(Y,0)].\]  

\smallbreak\n(ii) For $A,B\in\vg$, with a similar computation to the one above,
we obtain
\[ J[(0,A),(0,B)]=[J(0,A),(0,B)]+[(0,A),J(0,B)]+J[J(0,A),J(0,B)].\]  

\smallbreak\n(iii) Finally, for $X\in\ug,\,A\in\vg$, one obtains
\begin{align*}
J[(X,0),(0,A)] & = J(-\mu(A)X,\rho(X)A)\\
               & = (-\varphi^{-1}(\rho(X)A),-\varphi(\mu(A)X))\\
               & = (-X\cdot\varphi^{-1}A,-A\cdot\varphi X),\\
[J(X,0),(0,A)] & = [(0,\varphi X),(0,A)] = (0,[\varphi X,A]), \\
[(X,0),J(0,A)] & = [(X,0),(-\varphi^{-1}A,0)]=(-[X,\varphi^{-1}A],0),\\
\end{align*}
and also
\begin{align*}
J[J(X,0),J(0,A)] & = J[(0,\varphi X),(-\varphi^{-1}A,0)] \\
                & = -J[(-\varphi^{-1}A,0),(0,\varphi X)]\\
     & = -J(\mu(\varphi X)\varphi^{-1}A,-\rho(\varphi^{-1}A)\varphi X)\\
     & = -J(\varphi^{-1}(\varphi X\cdot A),-\varphi(\varphi^{-1}A\cdot X))\\
                & = (-\varphi^{-1}A\cdot X,-\varphi X\cdot A).
\end{align*}
Using (\ref{tfree}) again, we obtain 
\[ J[(X,0),(0,A)]=[J(X,0),(0,A)]+[(X,0),J(0,A)]+J[J(X,0),J(0,A)].\]  
This concludes the verification of the integrability of $J$. Thus, $J$ is a
complex structure which anticommutes with $E$, since $\vg=J\ug$, showing that
$\{J,\,E\}$ is a complex product structure on $\ggo$.
\end{proof}

\begin{rems} (i) Suppose that the LSA structure on $\vg$ is trivial, i.e.\
$A\cdot B=0$ for all $A,B\in\vg$.  Hence $\vg$ is an abelian ideal of
$\ug\bowtie\vg$ and, in fact, $\ug\bowtie\vg\cong\aff(\ug)$.

\smallbreak\n(ii) Majid \cite{Maj} has constructed for the Lie algebra $\ggo$
of a compact semisimple Lie group another Lie algebra $\ggo'$ such that
$(\ggo,\ggo')$ determines a matched pair of Lie algebras.  We have shown that
for each Lie algebra carrying an LSA structure there exists another Lie algebra
of the same dimension (namely, the abelian one) such that they determine a
matched pair of Lie algebras.
\end{rems}

\vs

\begin{cor}\label{Rn}
Consider the abelian Lie algebra $\RR^n$ endowed with an LSA structure. Then
$(\RR^n,\RR^n,\rho,\rho)$ is a matched pair of Lie algebras, where
$\rho:\RR^n\longrightarrow\glg(\RR^n)$ is the representation of $\RR^n$ given
by $\rho(X)Y=X\cdot Y$. Hence, $\ggo=\RR^n\bowtie\RR^n$ admits a complex
product structure.
\end{cor}

\begin{proof}
Simply set $\varphi=\Id_{\RR^n}$ in Theorem \ref{converse}. It is trivial to
verify that (\ref{jacobi1}) and (\ref{jacobi2}) hold and then the corollary
follows.
\end{proof}

The previous corollary can be generalized in the following way.

\vs

\begin{cor}
Suppose the abelian Lie algebra $\RR^n$ carries {\em two} LSA structures,
denoted by $\RR^n\times\RR^n\longrightarrow\RR^n,\,(X,Y)\mapsto X\cdot Y$ and
$\RR^n\times\RR^n\longrightarrow\RR^n,\,(X,Y)\mapsto X*Y$. Denote by
$L:\RR^n\longrightarrow\glg(\RR^n)$ and $L^*:\RR^n\longrightarrow \glg(\RR^n)$
the representations given by left multiplications: $L_XY=X\cdot Y,\quad
L^*_XY=X*Y$. If
\[ L_XL^*_Y=L_YL^*_X,\quad L^*_XL_Y=L^*_YL_X \]
for all $X,Y\in\RR^n$, then $(\RR^n,\RR^n,L,L^*)$ is a matched pair of Lie
algebras and hence, $\ggo=\RR^n\bowtie\RR^n$ carries a complex product
structure.
\end{cor}

\begin{proof}
Take $\varphi=\Id_{\RR^n}$ in Theorem \ref{converse}.
\end{proof}

Let $(\ggo,J,E)$ be a Lie algebra with a complex product structure and let
$(\ggo,\ggo_+,\ggo_-)$ be its associated double Lie algebra. We have shown that
both $\ggo_+$ and $\ggo_-$ carry LSA structures. A natural question to ask is
whether $\ggo$ itself supports an LSA structure, which restricts to the ones
already existing on $\ggo_{\pm}$. We can naturally extend the bilinear products
on $\ggo_+$ and $\ggo_-$ to $\ggo$ by using the representations
$\rho,\,\mu,\,\tilde{\rho},\,\tilde{\mu}$ in the following fashion:
\[ (X+A)\cdot(Y+B)=\tilde{\rho}(X)Y+\rho(X)B+\mu(A)Y+\tilde{\mu}(A)B,\quad 
X,Y\in\ggo_+,A,B\in\ggo_-, \] or more simply
\begin{equation} \label{noflat} 
(X+A)\cdot(Y+B)=X\cdot Y+\rho(X)B+\mu(A)Y+A\cdot B 
\end{equation}
where we are using the LSA structure on $\ggo_{\pm}$. It is a simple matter to
verify that this product on $\ggo$ verifies (\ref{tfree}). However, it does not
necessarily fulfill (\ref{flat}). We will give an example of this fact in the
next section (see Example \ref{nonflat}). 

To see when (\ref{noflat}) satisfies (\ref{flat}), we introduce the following
notation. Define $\Phi:\ggo_+\longrightarrow
(\ggo_-\otimes\ggo_-)^*\otimes\ggo_-$ by
\begin{equation}\label{phi}
\Phi(X)(A,B)=\rho(X)(A\cdot B)-(\rho(X)A)\cdot
B-A\cdot(\rho(X)B)+\rho(\mu(A)X)B
\end{equation}
and $\Psi:\ggo_-\longrightarrow(\ggo_+\otimes\ggo_+)^*\otimes\ggo_+$ by
\begin{equation}\label{psi}
\Psi(A)(X,Y)=\mu(A)(X\cdot Y)-(\mu(A)X)\cdot Y-X\cdot(\mu(A)Y)+\mu(\rho(X)A)Y,
\end{equation}
for $X,Y\in\ggo_+,\,A,B\in\ggo_-$. Then equations (\ref{jacobi1}) and
(\ref{jacobi2}) can be rewritten simply as
\[ \Phi(X)(A,B)=\Phi(X)(B,A),\quad \Psi(A)(X,Y)=\Psi(A)(Y,X),\]
for $X,Y\in\ggo_+,\,A,B\in\ggo_-$. Using Lemma 3.0.3 of \cite{DM}, we obtain
that the product defined above determines an LSA structure on $\ggo$ if and
only if $\Phi\equiv 0$ and $\Psi\equiv 0$.

Summing up,

\vs

\begin{prop} 
Let $(\ggo,J,E)$ be a Lie algebra with a complex product structure with
associated double Lie algebra $(\ggo,\ggo_+,\ggo_-)$ and
$\ggo_-=J\ggo_+$. There exists an LSA structure on $\ggo$ extending the LSA
structures on $\ggo_+$ and $\ggo_-$ if and only if $\Phi\equiv 0$ and
$\Psi\equiv 0$, where $\Phi$ and $\Psi$ are as in (\ref{phi}) and (\ref{psi}).
\end{prop}

\vs

\begin{cor}
Let $(\ggo,J,E)$ be a Lie algebra with a complex product structure and let
$(\ggo,\ggo_+,\ggo_-)$ be its associated double Lie algebra, with
$\ggo_-=J\ggo_+$. Suppose that $\ggo_-$ is an ideal of $\ggo$. Then $\ggo$
admits an LSA structure which extends the LSA structures on $\ggo_+$ and
$\ggo_-$
\end{cor}

\begin{proof} We already know that in this case the LSA structure on $\ggo_-$
is trivial (i.e.\ $A\cdot B=0$ for all $A,B\in\ggo_-$) and consequently, the
$\mu$ is also trivial.  Substituting into (\ref{phi}) and (\ref{psi}), we get
$\Phi\equiv 0$ and $\Psi\equiv 0$. Hence $\ggo$ does admit an LSA structure.
\end{proof}

\vs

\section{Torsion-free connections and $G$-structures}

We start recalling that for an arbitrary connection $\nabla$ on (the tangent
bundle of) a manifold $M$, the torsion and curvature tensor fields $T$ and $R$
are defined by \begin{gather*} T(X,Y)=\nabla_XY-\nabla_YX-[X,Y]\\R(X,Y)=
\nabla_X\nabla_Y-\nabla_Y\nabla_X-\nabla_{[X,Y]} \end{gather*} for $X,Y$
smooth vector fields on $M$. The connection is called {\em torsion-free} when
$T=0$, and {\em flat} when $R=0$. A flat torsion-free affine connection on $M$
gives rise to an {\em affine structure} on $M$ \cite{AUM,Mil}.

Let $G$ be a Lie group with Lie algebra $\ggo$ and suppose that $G$ admits a
left-invariant connection $\nabla$. This means that if $X,Y\in\ggo$ are two
left-invariant vector fields on $G$ then $\nabla_XY\in\ggo$ is also
left-invariant. Accordingly, one may define a connection on a Lie algebra
$\ggo$ to be merely a $\ggo$-valued bilinear form $\ggo\times\ggo
\longrightarrow\ggo$. One can speak of the torsion and curvature of such a
connection using the formulae above, with brackets determined by the structure
of $\ggo$. We can now re-interpret the results obtained in previous sections
in terms of connections.

\smallbreak

Let $\ggo$ be a Lie algebra with a complex product structure $\{J,E\}$ and let
$(\ggo,\ggo_+,\ggo_-)$ be the associated double Lie algebra. In \S4, we have
shown that both $\ggo_+$ and $\ggo_-$ carry an LSA structure, and we
introduced a bilinear product on $\ggo$ which extends these LSA structures and
satisfies (\ref{tfree}). This product becomes a torsion-free affine connection
$\ncp$ (CP for `complex product') on $\ggo$ by setting \[\ncp_XY=X\cdot Y,\]
for $X,Y\in \ggo$. The torsion-free condition reads \[X\cdot Y-Y\cdot
X=[X,Y],\] and the connection is flat if \[X\cdot(Y\cdot Z)-Y\cdot(X\cdot
Z)=[X,Y]\cdot Z.\] Since this is exactly condition (\ref{flat}), it follows
that $\ncp$ restricts to {\em flat} torsion-free connections on $\ggo_+$ and
$\ggo_-$. In general though, $\ncp$ is not itself flat (see
Example~\ref{nonflat} in \S6).

\vs

\begin{prop}\label{unique} If $(\ggo,J,E)$ is a Lie algebra carrying a complex
product structure then $\ncp J=0=\ncp E$. Moreover, $\ncp$ is the unique
torsion-free connection on $\ggo$ for which $J$ and $E$ are parallel.
\end{prop}

\begin{proof} The parallelism of $J,E$ is quickly established from the 
definitions of the representations $\rho,\,\mu,\,\tilde{\rho},\,\tilde{\mu}$
and equation (\ref{noflat}).

As to uniqueness, suppose $\nabla,\nabla'$ are two connections on $\ggo$ which
satisfy the conditions in the statement. Define $A:\ggo\times\ggo
\longrightarrow\ggo$ by $A(X,Y)=\nabla_XY-\nabla'_XY$. Since $\nabla,\nabla'$
are both torsion-free, \[\begin{array}{rcl} A(Y,X) &=&\nabla_YX-\nabla'_YX\\
&=&\nabla_XY+[Y,X]-\nabla'_XY+[X,Y]\ =\ A(X,Y),\end{array}\] and $A$ is
symmetric. Since $\nabla J=0$, we obtain \[\begin{array}{rcl} A(X,JY)&=&
\nabla_XJY-\nabla'_XJY\\&=& J(\nabla_XY-\nabla'_XY)\ =\ JA(X,Y)\end{array}\]
and (by symmetry) $A(JX,Y)=JA(X,Y)$. Since also $\nabla E=0$, the following is
valid: \[A(X,EY)= EA(X,Y)=A(EX,Y).\] Then, $A(JX,EY)=EA(JX,Y)=JEA(X,Y)$, and
\[A(JX,EY)=JA(X,EY)=EJA(X,Y)=-JEA(X,Y).\] It follows that $A(X,Y)=0$
for all $X,\,Y\in\ggo$, and so $\nabla=\nabla'$.\end{proof}

We have shown in Theorem \ref{hyper} that the complex product structure
$\{J,E\}$ on $\ggo$ determines a hypercomplex structure $\{\II,\JJ\}$ on
$\hat\ggo=(\ggo^\CC)_\RR$. In this case, any Lie group $G$ with $\Lie(G)=
\hat\ggo$ is a {\em hypercomplex manifold}, since it admits a pair
$\{\II,\JJ\}$ of anticommuting complex structures. These are obtained by
left-translating those that are already defined on $\ggo\cong T_eG$. Every
hypercomplex structure on $M$ uniquely determines a torsion-free connection
$\nhc$ (HC for `hypercomplex'), sometimes called the {\em Obata connection},
with respect to which $\II,\JJ$ are parallel. Thus, \[\nhc\II=0,\quad
\nhc\JJ=0\] (see \cite{Ob}). Since $\II,\JJ$ are left-invariant, it follows
that $\nhc$ is left-invariant. In this way, we obtain a torsion-free
connection on $\hat\ggo$, that we also call the Obata connection.

On the other hand, the connection $\ncp$ on $\ggo$ extends naturally to a
torsion-free connection $\wncp$ on $\hat\ggo$ by setting \[\wncp_{X+IY}
(X'+IY')=(\ncp_XX'-\ncp_YY')+I(\ncp_XY'+\ncp_YX'),\] for $X,X',Y,Y'\in\ggo$.
Since $\ncp J=0$, it is clear that $\wncp\JJ=0$, and a simple computation shows
that $\wncp\II=0$. By uniqueness of the Obata connection, we have
$\wncp=\nhc$. We have thus proven

\medskip

\begin{prop} The canonical extension of $\ncp$ to $\hat\ggo$ coincides
with the connection $\nhc$ determined by the hypercomplex structure
$\{\II,\JJ\}$ on $\hat\ggo$ given by Theorem \ref{hyper}.\end{prop}

\smallskip

\begin{cor}\label{obata} The Obata connection $\nhc$ on $\hat\ggo$ is flat
if and only if the complex product connection $\ncp$ on $\ggo$ is
flat.\end{cor}

\smallskip

\begin{exam} The corollary applies to $\ggo=\glg(2n,\RR)$. For the proof of
Proposition~\ref{gl} shows that the complex product structure on $\glg(2n,\RR)$
realizes it as $2n$ copies of the flat structure on $\RR^{2n}$. In the same
way, the induced hypercomplex structure on $\glg(2n,\CC)$ is $2n$ copies
(\ref{2n}) of flat space $\HH^n$.\end{exam}

\medskip

Recall the matrices $J_0,\,E_0$ defined in (\ref{E0J0}), and acting by left
multiplication on the space $\RR^{2n}$ of column vectors. The subgroup of
$GL(2n,\RR)$ consisting of matrices commuting with $J$ can be identified with
$GL(n,\CC)$. The subgroup of those commuting with both $J$ and $E$ consists of
invertible block-diagonal matrices \begin{equation}\label{block}
\left(\begin{array}{cc}A&0\\0&A\end{array}\right)\end{equation} and can be
identified with $GL(n,\RR)$.

We can generalize the preceding theory by

\medskip

\begin{defi} 
A \textit{complex product manifold} is a smooth manifold $M$, of even
dimension $2n$, equipped with a $GL(n,\RR)$-structure (relative to the above
inclusion $GL(n,\RR)\subset GL(2n,\RR)$) admitting a $GL(n,\RR)$ connection.
\end{defi}

In the general context of $G$-structures, there is a well-known series of
obstructions to a given structure being equivalent to the standard structure
on flat space \cite{Stern}. The existence of a torsion-free connection is in
general only the first such obstruction. The proof of Proposition~\ref{unique}
can be extended to show that a torsion-free $GL(n,\RR)$-connection will (if it
exists) be unique. Full details of this and the subsequent remarks appear in
\cite{A}.

Given a $GL(n,\RR)$-structure, we can regard its Lie algebra $\glg(n,\RR)$ as
a subspace of $\End T=T\otimes T^*$, where $T=\RR^{2n}$ stands for the tangent
space representation. Torsion properties are determined by the natural mapping
\[f_1:\glg(n,\RR)\otimes T^*\to T\otimes\alt^2T^*.\] The difference of two
compatible torsion-free connections is a tensor that hypothetically takes
values at each point in $\ker f_1$, which is the so-called the first
prolongation of $\glg(n,\RR)$ (as a subalgebra of $\glg(2n,\RR)$). The above
remarks then imply

\vs

\begin{lem}
$\ker f_1=\{0\}$.
\end{lem}

\vs

Let $T=T_+\oplus T_-$ be the decomposition reflecting (\ref{block}), so that
$T_+\cong T_-$ is the standard representation of $GL(n,\RR)$. The torsion-free
condition ensures that the distributions determined by $T_+,\,T_-$ are
integrable and totally geodesic. Moreover, $M$ has an integrable complex
structure $J$ for which $T_-=JT_+$. The spaces $T_\pm$ now play the role of
$\ggo_\pm$ in \S4.

The second (and only remaining) obstruction to flatness of the
$GL(n,\RR)$-structure of $M$ is provided by the curvature tensor $R$ of
$\nabla$. This takes values in the kernel of the linear mapping
\[f_2:\glg(n,\RR)\otimes \alt^2T^*\to T\otimes\alt^3T^*\] that establishes the
first Bianchi identity. The existence of an LSA structure on $\ggo_\pm$ can be
seen as a special case of the next result, that can be proved by working in
terms of the standard $GL(n,\RR)$-module $V=\RR^{n}$, using the isomorphisms
$T_+\cong V\cong T_-$ and $\glg(n,\RR)\cong V^*\otimes V$ \cite{A}.

\vs

\begin{lem} The restriction of $R$ to $\glg(n,\RR)\otimes\alt^2T^*_\pm$
vanishes.
\end{lem}

\smallbreak\n This implies that the restriction of $\nabla$ to each of the
integrable submanifolds tangent to $T_\pm$ is flat.

The symmetries of $R$ resemble those of the curvature tensor of a hypercomplex
manifold (described in \cite{Sal}), which is not surprising given the intimate
relationship between the two structures. Hypercomplex structures are
associated to the group $GL(n,\HH)$ in the same way that complex product
structures are associated to $GL(n,\RR)$. The two groups can be extended to
$GL(n,\RR)\times GL(2,\RR)$ and $GL(n,\HH)\times GL(1,\HH)$ respectively, and
the resulting $G$-structures are both types of paraconformal geometries
\cite{BE}.

\smallbreak

The above manifold theory can of course be applied to any Lie group $G$ whose
Lie algebra $\ggo$ admits a complex product structure $\{J,E\}$. We may regard
$J,E$ as tensors on $G$ by left-translating the ones on $\ggo=T_eG$. If
$(\ggo,\ggo_+,\ggo_-)$ is the double Lie algebra defined by $E$, let $G_+$ and
$G_-$ be Lie subgroups of $G$ such that $\Lie(G_\pm)=\ggo_\pm$. We learn from
\cite{LW} that $(G,G_+,G_-)$ is a local double Lie group, i.e.\ the mapping
$G_+\times G_-\longrightarrow G$ defined by $(g_+,g_-)\mapsto g_+g_-$ is a
local diffeomorphism near the identities. Since $\ggo_\pm$ carry LSA
structures, both $G_+$ and $G_-$ are flat affine Lie groups, whilst $G$ itself
need not be flat.

In the notation of \S2, let $\hat G$ be a Lie group with Lie algebra
$\hat\ggo$. Then $\hat G$ acquires a $GL(n,\CC)$-structure, obtained by
complexifying (\ref{block}) and then embedding the complex matrix in
$GL(4n,\RR)$. The induced hypercomplex structure on $\hat G$ comes about as a
result of inclusions \[GL(n,\CC)\subset GL(n,\HH)\subset GL(4n,\RR).\]

\smallbreak Before concluding the $G$-structure approach, we consider cases in
which a complex product structure reduces from $GL(n,\RR)$ to either (i) the
orthogonal group $O(n)$, or (ii) the symplectic group $Sp(n/2,\RR)$.

\smallbreak\n(i) $O(n)$ represents the subgroup preserving a scalar product
$g$ on $T_+$ or $T_-$. These two representations are equivalent by means of
$J$, so we extend $g$ to $T=\RR^{2n}$ by the requirement that
\[g(JX,JY)=g(X,Y),\qquad g(EX,EY)=g(X,Y),\] for all tangent vectors $X,Y$. A
2-form is now defined in the usual way, by setting \[\omega(X,Y)=g(X,JY).\] If
this is closed on a complex product manifold, then $g$ is a reducible {\em
K\"ahler metric}, and $\ncp$ coincides with the Levi-Civita connection
$\nabla^g$. It is however possible to drop the integrability assumption on $J$,
and consider the complementary {\em Lagrangian foliations} defined by the
distributions $T_\pm$ on a symplectic manifold.

\smallbreak\n(ii) $Sp(n/2,\RR)$ represents the group preserving a
non-degenerate 2-form $\omega_1$ on $T_+$ or $T_-$, each of which must then
have even dimension $n=2m$. We may extend $\omega_1$ to $T$ by the requirement
that \[\omega_1(JX,JY)=\omega(X,Y),\qquad \omega_1(EX,EY)= \omega_1(X,Y),\] in
analogy to $g$ above. A symmetric bilinear form $h$ is now defined by setting
\[ h(X,Y)=\omega_1(JX,Y),\] and it follows that $h(EX,EY)=-h(X,Y)$; thus
$T_\pm$ are isotropic and $h$ has signature $(2m,2m)$. Additional 2-forms are
defined by setting \[\omega_2(X,Y)=h(X,EY),\qquad \omega_3(X,Y)= h(X,JEY).\] 

\smallbreak A manifold is called {\em hypersymplectic} if it has such an
$Sp(n/2,\RR)$ structure for which all three forms $\omega_i$ are closed
\cite{Hi}. Such manifolds are also called {\em neutral hyperk\"ahler}
\cite{Kam,FPPS}, and admit tensors $J,E,F$ satisfying (\ref{rules}). The
complex 2-form $\Omega=\omega_2+i\omega_3$ has type $(2,0)$ relative to $J$
and the fact that $\Omega^m$ is closed and non-zero implies that $J$ is
integrable. It follows that $h$ is pseudo-K\"ahler, and that $\ncp$ coincides
with $\nabla^h$. The subspaces $T_\pm$ are Lagrangian relative to $\omega_2$,
though it is more usual to seek submanifolds that are special Lagrangian
relative to the forms $\omega_1$ and $\hbox{Im}\,(\Omega^m)$ (meaning that these
restrict to zero).

In conclusion, a hypersymplectic manifold always has an underlying complex
product structure (with $n$ even), in the same sort of way that a
hyperk\"ahler manifold has a subordinate hypercomplex structure.

\vs

\section{Complex product structures on 4-dimensional Lie algebras}
\label{examples}

In this section we shall consider some concrete examples of 4-dimensional Lie
algebras carrying complex product structures. We will determine the associated
isomorphisms with bicrossproduct Lie algebras and also the corresponding local
double Lie groups.

Let us first make the following observation. If a 4-dimensional Lie algebra
admits a complex product structure, then it can be written as the sum of two
2-dimensional subalgebras. We will use the fact that any 2-dimensional Lie
algebra is either abelian or isomorphic to $\aff(\RR)$ (see Example
\ref{ejemplo} $\rii$).

We will consider separately the solvable Lie algebras and the non-solvable Lie
algebras. It is easy to see that there are only two of
the latter, using the Levi decomposition and the classification of 3-dimensional
simple Lie algebras. The only possibilities are the reductive Lie algebras $\RR
\oplus \sog(3)$ and $\RR \oplus \slg(2,\RR)\cong \glg(2,\RR)$.

\subsection{Reductive case}

\begin{exam} Let us begin with $\RR\oplus\slg(2,\RR)\cong \glg(2,\RR)$. This
Lie algebra has a basis $\{W,X,Y,Z\}$ with $Z$ central and non-zero brackets
given by \[ [W,X]=2X, \quad [W,Y]=-2Y, \quad [X,Y]=W. \] We already know from
Proposition \ref{gl} that $\glg(2,\RR)$ carries a complex product structure
$\{J,E\}$ given by right multiplication with the standard almost complex and
almost product structures $J_0$ and $E_0$ on $\RR^{2n}$. Identifying $W,X,Y,Z$
with $(2\times 2)$ matrices, it is easy to see that \[JW=-(X+Y),\quad
JX=\ft(W+Z),\quad JY=\ft(W-Z),\quad JZ=-X+Y,\] while
\[ EW=Z, \quad EX=-X, \quad EY=Y, \quad EZ=W. \] The eigenspaces for $E$ are
the following subalgebras of $\glg(2,\RR)$: \[
\ggo_+=\text{span}\{Y,W+Z\},\quad \ggo_-=\text{span}\{X,W-Z\}.\] Observe that
both $\ggo_+$ and $\ggo_-$ are isomorphic to $\glg(2,\RR)$; therefore the
associated double Lie algebra is $(\glg(2,\RR,\aff(\RR),\aff(\RR))$ and its
local double Lie group is $(GL(2,\RR),G_+,G_-)$, where both $G_+$ and $G_-$
are locally isomorphic to $\Aff(\RR)$.  \end{exam}

\smallskip

\begin{rem} Sasaki classified the complex structures on $\glg(2,\RR)$ in
\cite{Sa}. They are parametrized by $d\in\CC\setminus\{0\}$ on the curve
$\operatorname{Re}(1/d)=-1$. The complex structure which appears above
corresponds to the case $d=-1$.  \end{rem}

Now consider the other $4$-dimensional reductive Lie algebra.

\smallskip 

\begin{prop}
The Lie algebra $\RR \oplus \sog(3)$ does not support any complex product
structure.
\end{prop}

\begin{proof} The compact Lie algebra $\ggo=\RR\oplus \sog(3)$ admits an
invariant inner product $\langle\,,\rangle$. Consider a 2-dimensional Lie
subalgebra $\ug$ of $\ggo$. As $\ug$ is a 2-dimensional Lie algebra, we know
that $\ug$ is abelian or isomorphic to $\aff(\RR)$. If $\ug\cong\aff(\RR)$,
there is a basis $\{U_1,U_2\}$ of $\ug$ with $[U_1,U_2]=U_2$. Thus \[ \langle
[U_1,U_2],U_2\rangle = \langle U_2,U_2\rangle=||U_2||^2,\] but, using the
invariance of $\langle\,,\rangle$, \[ \langle [U_1,U_2],U_2\rangle = -\langle
U_2,[U_1,U_2]\rangle,\] from where we obtain $\langle
[U_1,U_2],U_2\rangle=0$. Hence, $||U_2||^2=0$, a contradiction. Therefore,
$\ug$ must be abelian.

Let us now suppose that $\ggo$ admits a complex product structure $\{J,E\}$
with associated double Lie algebra $(\ggo,\ggo_+,\ggo_-)$. As we have just
seen, both $\ggo_+$ and $\ggo_-$ are abelian.  For $U,U'\in\ggo_+$ and
$V,V'\in\ggo_-$, we compute
\[ \langle [U,V],U'\rangle=-\langle V,[U,U']\rangle = 0 \]
and 
\[ \langle [U,V],V'\rangle=-\langle [V,U],V'\rangle=-\langle U,[V,V']\rangle = 0 \]
because $\ggo_{\pm}$ are abelian. Hence $[U,V]=0$ and then $\ggo$ should be
abelian, which is false. The proof is now complete.
\end{proof}

\subsection{Solvable case}

\begin{exam}\label{nonflat} Let $\ggo$ be the Lie algebra defined by
$\ggo=\text{span}\{A,B,C,D\}$ with non-zero Lie bracket relations given by \[
[A,B]=B,\,[A,C]=-C,\,[A,D]=-D.\] $\ggo$ lies in the class $A2$ of the
classification made in \cite{O}, with $\lambda=-1$. It follows that $\ggo$
admits only one complex structure $J$, up to isomorphism, given by \[ JA=B,\,
JC=D.\] Consider the following subalgebras of $\ggo$: \[
\ggo_+=\text{span}\{A-D, C\},\quad \ggo_-=\text{span}\{B+C, D\}. \] Set
$E|_{\ggo_+}=\Id,\;E|_{\ggo_-}=-\Id$; the endomorphism $E$ of $\ggo$ is clearly a
product structure on $\ggo$ which anticommutes with the complex structure $J$,
giving rise then to a complex product structure on $\ggo$. Thus, there exists
an LSA structure on $\ggo_+$ and $\ggo_-$, which can be extended to a bilinear
product on $\ggo$. We already know that this product satisfies condition
(\ref{tfree}) but we will show that it does not satisfy (\ref{flat}). Denote
by $L_v$ the endomorphism of $\ggo$ given by left-multiplication with
$v\in\ggo$. Using equation (\ref{noflat}), we obtain \[ \begin{array}{ll}
L_{A-D}(A-D)=A-D,\qquad & L_{B+C}(A-D)=2C, \\ L_{A-D}(C)=-C, & L_{B+C}(C)=0, \\
L_{A-D}(B+C)=B+C, & L_{B+C}(B+C)=2D, \\ L_{A-D}(D)=-D, & L_{B+C}(D)=0. \\
\end{array} \] and $L_C\equiv 0,\;L_D\equiv 0$. Now, setting
$x=A-D,\,y=B+C,\,z=A-D$ in (\ref{flat}), we get \[ x\cdot(y\cdot z)-(x\cdot
y)\cdot z=-4C, \] while \[ y\cdot(x\cdot z)-(y\cdot x)\cdot z=2C, \] and hence
(\ref{flat}) does not hold.  \end{exam}

\vs 

\begin{exam}\label{S1}
Consider the Lie algebra defined by $\ggo=\text{span}\{X,Y,Z,W\}$ with the only
non-zero bracket given by $[X,Y]=Z$. This Lie algebra is isomorphic to
$\hg_3\oplus\RR$, a direct sum of ideals, where $\hg_3$ denotes the
$3$-dimensional Heisenberg Lie algebra. This Lie algebra belongs to the class
$S1$ in the classification made in \cite{Sn}, from where it follows that there
is, up to equivalence, only one complex structure on $\ggo$ and it is given by
\[ JX=Y,\quad JZ=W. \]

\vs

\begin{prop} 
$\hg_3\oplus\RR$ admits complex product structures. With each of these
structures, the associated double Lie algebra is
$(\hg_3\oplus\RR,\RR^2,\RR^2)$ and the associated local double Lie group is
$(G,G_+,G_-)$, where $G$ is locally isomorphic to $H_3\times\RR$ and $G_+,G_-$
are locally isomorphic to $\RR^2$.  
\end{prop}

\begin{proof} It can be shown that the product structures on $\ggo$ which 
anticommute with $J$ are parametrized by
\[ E_{\theta}=\begin{pmatrix}
                      1 & 0 & 0 & 0 \cr 0 & -1 & 0 & 0 \cr 0 & 0 & \cos\theta &
                      \sin\theta \cr 0 & 0 & \sin\theta & -\cos\theta \cr
                      \end{pmatrix} ,\] 
where the ordered basis used for the matrix representation is
$\{X,Y,Z,W\}$. The pairs $\{J,E_{\theta}\}$ with $\theta\in[0,2\pi)$ exhaust
all the complex product structures on $\ggo$, up to equivalence.

For a fixed $\theta\in[0,2\pi)$, let us establish the decomposition of $\ggo$
into a sum of subalgebras. In order to do so, we simply have to determine the
eigenspaces of $E_{\theta}$. The eigenspace $\ggo_+$ corresponding to the
eigenvalue $+1$ is given by
\[ \ggo_+=\text{span}\{X,\,\ctt Z+\stt W\},\]
while the eigenspace $\ggo_-$ corresponding to the eigenspace $-1$ is given by
\[ \ggo_-=\text{span}\{Y,\,-\stt Z+\ctt W\}.\]
(We indicate the trigonometric arguments as subscripts for visual clarity.)
Observe that both $\ggo_+$ and $\ggo_-$ are abelian Lie algebras. Hence, we
have $\ggo=\RR^2\bowtie\RR^2$ for suitable representations and certain LSA
structures on each $\RR^2$.
\end{proof}
\end{exam}

\smallskip

\begin{exam}\label{S3} Let $\ggo=\text{span}\{A,B,C,D\}$ be the Lie algebra
with non-zero Lie bracket relations given by \[[A,B]=B,\quad[A,C]=C,\quad
[A,D]=D.\] Then $\ggo$ lies in the class $A4$ of the classification \cite{O},
and is the Lie algebra corresponding to the real hyperbolic space $\RR H^4$,
i.e.\ the simply connected Lie group $S$ with $\Lie(S)=\ggo$ acts simply
transitively on $\RR H^4$. It follows that $\ggo$ admits only one complex
structure $J$, up to isomorphism, given by \[ JA=B, \quad JC=D.\]

\smallskip

\begin{prop} 
$\ggo$ admits complex product structures. The possible associated double Lie
algebras are either $(\ggo,\aff(\RR),\aff(\RR))$ or $(\ggo,\RR^2,\aff(\RR))$
and the associated local double Lie groups are $(S,G_+,G_-)$, where both
$G_+,G_-$ are both locally isomorphic to $\Aff(\RR)$, or $G_+$ is locally
isomorphic to $\RR^2$ and $G_-$ to $\Aff(\RR)$.
\end{prop}

\begin{proof}
The following are examples of product structures on $\ggo$ which anticommute
with the complex structure $J$:
\[ E_{\theta}=\begin{pmatrix}
                      \cos\theta & \sin\theta & 0 & 0 \cr
                      \sin\theta & -\cos\theta & 0 & 0 \cr
                      0 & 0 & 1 & 0 \cr 
                      0 & 0 & 0 & -1 \cr
                     \end{pmatrix} ,\qquad 
   E'_{\theta}=\begin{pmatrix}
                      \cos\theta & \sin\theta & 0 & 0 \cr
                      \sin\theta & -\cos\theta & 0 & 0 \cr
                      -\sin\theta & 1+\cos\theta & 1 & 0 \cr
                      1+\cos\theta & \sin\theta & 0 & -1 \cr
                     \end{pmatrix} ,\]
\[ E''_{\theta}=\begin{pmatrix}
                      \cos\theta & \sin\theta & -\sin\theta & 1+\cos\theta \cr
                      \sin\theta & -\cos\theta & 1+\cos\theta & \sin\theta \cr
                      0 & 0 & 1 & 0 \cr
                      0 & 0 & 0 & -1 \cr
                     \end{pmatrix} ,\quad
   \tilde{E}=\begin{pmatrix}
                        -1 & 0 & 0 & 0 \cr 0 & 1 & 0 & 0 \cr -2 & 0 & 1 & 0 \cr
                        0 & 2 & 0 & -1 \cr \end{pmatrix}. \] 
The ordered basis used for the matrix representation is $\{A,B,C,D\}$. It can
be shown that every complex product structure on $\ggo$ is equivalent to
$\{J,E\}$ where $J$ is as above and $E\in \{E_{\theta}:\theta\in[0,2\pi)\}
\cup\{E'_{\theta}:\theta\in[0,2\pi),\theta\neq\pi\}
\cup\{E''_{\theta}:\theta\in[0,2\pi),\theta\neq\pi\}\cup\{\tilde{E}\}$,
and these are all inequivalent.

For $\theta\in[0,2\pi)$, let us compute the eigenspaces corresponding to
$E_{\theta}$. The eigenspace $\ggo_+$ corresponding to the eigenvalue $+1$ is
given by
\[ \ggo_+=\text{span}\{C,\,\ctt A+\stt B\},\]
while the eigenspace $\ggo_-$ corresponding to the eigenvalue $-1$ is
\[ \ggo_-=\text{span}\{D,\,-\stt A+\ctt B\}.\]
Observe that both $\ggo_+$ and $\ggo_-$ are isomorphic to $\aff(\RR)$, hence we
have an isomorphism $\ggo\cong\aff(\RR)\bowtie\aff(\RR)$, for suitable
representations and certain LSA structures on each $\aff(\RR)$.

For $\theta\in [0,2\pi),\,\theta\neq\pi$, let us consider the complex product
structures $\{J,E_{\theta}'\}$. The eigenspaces of $E_{\theta}'$ are
\[\begin{array}{c}
\ggo_+'=\text{span}\{C,\,\ctt A+\stt B+\ctt D\},\\[3pt]
\ggo_-'=\text{span}\{D,\,-\stt A+\ctt B-\ctt C\}.
\end{array}\]
Observe that both $\ggo_+$ and $\ggo_-$ are isomorphic to $\aff(\RR)$, hence 
$\ggo\cong\aff(\RR)\bowtie\aff(\RR)$, for suitable
representations and certain LSA structures on each $\aff(\RR)$.

For $\theta\in [0,2\pi),\,\theta\neq\pi$, let us consider the complex product
structures $\{J,E_{\theta}''\}$. The eigenspaces of $E_{\theta}''$ are
\[\begin{array}{c}
\ggo_+''=\text{span}\{\ctt A+\stt B,\,-\ctt A+
\stt C\},\\[3pt]
\ggo_-''=\text{span}\{-\stt A+\ctt B,\,-\ctt B+\stt D\}.
\end{array}\] 
Note that both $\ggo_+$ and $\ggo_-$ are isomorphic to $\aff(\RR)$, hence 
$\ggo\cong\aff(\RR)\bowtie\aff(\RR)$, for suitable
representations and certain LSA structures on each $\aff(\RR)$.

Finally, the eigenspaces for $\tilde{E}$ are \[
\tilde{\ggo}_+=\text{span}\{C,B+D\},\quad
\tilde{\ggo}_-=\text{span}\{D,A+C\}.\] Note that $\tilde{\ggo}_+$ is abelian
while $\tilde{\ggo}_-$ is isomorphic to $\aff(\RR)$. Thus, in this case we
have an isomorphism $\ggo\cong\RR^2\bowtie\aff(\RR)$, for suitable
representations and certain LSA structures on each summand.\end{proof}
\end{exam}

\vs

In contrast to above, not every solvable Lie algebra (among those carrying a
complex structure) admits a complex product structure.

\vs

\begin{exam}
Let $\ggo$ be the Lie algebra given by $\ggo=\text{span}\{A,X,Y,Z\}$ with Lie
bracket relations \[ [X,Y]=Z,\quad[A,X]=-Y,\quad[A,Y]=X. \]
This Lie algebra belongs to the class $H2$ of the classification made in
\cite{O} and it follows that it admits, up to equivalence, only one complex
structure $J$, which is given by \[ JA=Z,\quad JX=Y.\]
We will show that this Lie algebra does not admit any complex product
structure, using

\vs

\begin{lem}
$\ggo$ does not contain a non-abelian 2-dimensional Lie subalgebra.
\end{lem}

\begin{proof}
Suppose that $\ug$ is a non-abelian 2-dimensional Lie subalgebra of
$\ggo$. Then it is isomorphic to $\aff(\RR)$ and hence there is a basis
$\{U,V\}$ of $\ug$ such that $[U,V]=V$. If
$U=a_1A+x_1X+y_1Y+z_1Z,\,V=a_2A+x_2X+y_2Y+z_2Z$, the relation $[U,V]=V$ implies
at once that $a_2=0$ and
\[
\begin{cases}
a_1y_2=x_2,\\
a_1x_2=-y_2,\\
x_1y_2-y_1x_2=z_2.
\end{cases}
\]
From the first two equations we obtain that $y_2(a_1^2+1)=0$. Since we are
working over $\RR$, we have that $y_2=0$. But this implies that $x_2=0$ and
hence $z_2=0$, showing that $U=0$, which is a contradiction.
\end{proof}

Let us suppose now that $\ggo$ admits a complex product structure. Then, we
have $\ggo=\ggo_+\oplus\ggo_-$ with $\ggo_+$ and $\ggo_-=J\ggo_+$ subalgebras
of $\ggo$. From the lemma, we obtain that both $\ggo_{\pm}$ are abelian.  It
can be seen that in this case the complex structure $J$ satisfies the
condition \begin{equation}\label{nose} [JX,JY]=[X,Y] \end{equation} for all
$X,Y\in\ggo$.  A complex structure $J$ on a Lie algebra $\ggo$ satisfying
(\ref{nose}) is called {\em abelian}, and only solvable Lie algebras admit
such structures \cite{DF1}. A result in \cite{BD} states that if a Lie algebra
with commutator ideal of codimension 1 admits an abelian complex structure
then it is isomorphic to $\aff(\RR)$. Since $\ggo$ is not isomorphic to
$\aff(\RR)$, we obtain a contradiction and then $\ggo$ cannot carry a complex
product structure.  \end{exam}

\subsection{Induced hypercomplex structures} In this subsection we will give
some applications of Theorem~\ref{hyper}, using the examples studied in \S 6.1
and \S 6.2. We will use the following notation: for $X\in\ggo$, we denote
$\hat{X}=iX$ in $\ggo^{\CC}$.

\vs

\begin{exam}
Let us begin with $\glg(2,\RR)$. Consider on this Lie algebra the complex
product structure given in \S 6.1 and let $\{\II,\JJ\}$ be the hypercomplex
structure on $\glg(2,\CC)=\widehat{\glg(2,\RR)}$ given by
Theorem~\ref{hyper} and (\ref{2n}). $\glg(2,\CC)$ is an 8-dimensional real Lie algebra with a
basis $\{W,X,Y,Z,\hat{W},\hat{X},\hat{Y},\hat{Z}\}$ and Lie bracket given by
\begin{gather*}
[W,X]=-[\hat{W},\hat{X}]=2X,\quad [W,Y]=-[\hat{W},\hat{Y}]=-2Y,\quad [X,Y]=-
[\hat{X},\hat{Y}]=W, \\ [W,\hat{X}]=[\hat{W},X]=2\hat{X},\quad
[W,\hat{Y}]=-[\hat{W},Y]=-2\hat{Y},\quad [\hat{X},Y]=-[X,\hat{Y}]=\hat{W},
\end{gather*}
The hypercomplex structure $\{\II,\JJ\}$ is given by
\begin{gather*}
\II W=\hat{Z},\quad\II X=-\hat{X},\quad\II Y=\hat{Y},\quad\II Z=\hat{W},\\
\II \hat{W}=-Z,\quad\II \hat{X}=-X,\quad\II \hat{Y}=Y,\quad\II \hat{Z}=-W,
\end{gather*}
and
\begin{gather*}
\JJ W=-(X+Y),\quad \JJ X=\ft(W+Z),\quad \JJ Y=\ft(W-Z),\quad\JJ
Z=-X+Y,\\
\JJ \hat{W}=-(\hat{X}+\hat{Y}),\quad \JJ \hat{X}=\ft(\hat{W}+\hat{Z}),\quad \JJ
\hat{Y}=\ft(\hat{W}-\hat{Z}),\quad\JJ \hat{Z}=-\hat{X}+\hat{Y}.
\end{gather*}
\end{exam}

\vs 

\begin{exam}
We continue here Example \ref{nonflat}. $\hat{\ggo}=(\ggo^{\CC})_{\RR}$ is an
8-dimensional solvable Lie algebra with a basis
$\{A,B,C,D,\hat{A},\hat{B},\hat{C},\hat{D}\}$ and non-zero Lie brackets given
by
\begin{gather*}
[A,B]=-[\hat{A},\hat{B}]=B,\quad [A,C]=-[\hat{A},\hat{C}]=-C, \quad
[A,D]=-[\hat{A},\hat{D}]=-D, \\ [A,\hat{B}]=[\hat{A},B]=\hat{B}, \quad
[A,\hat{C}]=[\hat{A},C]=-\hat{C}, \quad [A,\hat{D}]=[\hat{A},D]=-\hat{D}.
\end{gather*}
The complex product structure considered on $\ggo$ gives rise to a hypercomplex
structure $\{\II,\JJ\}$ on $\hat{\ggo}$. This hypercomplex structure is given
by
\begin{gather*}
\II A=\hat{A}-2\hat{D},\quad \II B=-\hat{B}-2\hat{C},\quad \II C=\hat{C},\quad
\II D=-\hat{D},\\ \II\hat{A}=-A+2D,\quad \II \hat{B}=B+2C,\quad \II
\hat{C}=-C,\quad \II \hat{D}=D,
\end{gather*}
and 
\[ \JJ A=B,\quad \JJ C=D,\quad \JJ \hat{A}=\hat{B},\quad \JJ \hat{C}=\hat{D}.\]
According to Example \ref{nonflat} and Corollary \ref{obata} we obtain that the
Obata connection on $\hat{\ggo}$ (or on any Lie group $G$ with
$\Lie(G)=\hat{\ggo}$) is not flat.
\end{exam}

\vs

\begin{exam} We continue here Example \ref{S1}. If $\ggo=\hg_3\oplus\RR$, then 
$\hat{\ggo}=(\ggo^{\CC})_{\RR}$ is an $8$-dimensional 2-step nilpotent Lie
algebra with a basis $\{X,Y,Z,W,\hat{X},\hat{Y},\hat{Z},\hat{W}\}$ and non-zero
Lie brackets given by
\[ [X,Y]=-[\hat{X},\hat{Y}]=Z,\quad [X,\hat{Y}]=[\hat{X},Y]=\hat{Z}.\]
The hypercomplex structure $\{\II_{\theta},\JJ\}$ on $\hat{\ggo}$ given by 
Theorem~\ref{hyper} is 
\begin{gather*}
\II_{\theta}X=\hat{X}, \quad \II_{\theta}Y=-\hat{Y}, \quad
\II_{\theta}Z=\ct\hat{Z}+\st\hat{W}, \quad
\II_{\theta}W=\st\hat{Z}-\ct\hat{W},\\ \II_{\theta}\hat{X}=-X,
\quad \II_{\theta}\hat{Y}=Y, \quad \II_{\theta}\hat{Z}=-\ct Z-\st W, \quad 
\II_{\theta}\hat{W}=-\st Z+\ct W,\\
\JJ X=Y,\quad \JJ Z=W,\quad \JJ \hat{X}=\hat{Y},\quad \JJ \hat{Z}=\hat{W}.
\end{gather*}
\end{exam}

\vs

\begin{rem} A classification of 8-dimensional nilpotent Lie groups carrying an
abelian hypercomplex structure (meaning all its complex structures satisfy
(\ref{nose})) was given in \cite{DF}. According to this classification, if $G$
is a Lie group whose Lie algebra is the one considered in the previous
example, then $G$ is a isomorphic to a trivial extension of a group of type
$H$ with centre of dimension 2.  \end{rem}

\vs

\begin{exam} We continue here Example \ref{S3}. $\hat{\ggo}=(\ggo^{\CC})_{\RR}$ 
is an 8-dimensional solvable Lie algebra with a basis
$\{A,B,C,D,\hat{A},\hat{B},\hat{C},\hat{D}\}$ and non-zero Lie brackets given
by
\begin{gather*}
[A,B]=-[\hat{A},\hat{B}]=B,\quad [A,C]=-[\hat{A},\hat{C}]=C, \quad
[A,D]=-[\hat{A},\hat{D}]=D, \\ [A,\hat{B}]=[\hat{A},B]=\hat{B}, \quad
[A,\hat{C}]=[\hat{A},C]=\hat{C}, \quad [A,\hat{D}]=[\hat{A},D]=\hat{D}.
\end{gather*}
Each complex product structure on $\ggo$ gives rise to a hypercomplex structure
on $\hat{\ggo}$.

If we take $\{J,E_{\theta}\}$, we have the hypercomplex structure
$\{\II_{\theta},\,\JJ\}$, where
\begin{gather*}
\II_{\theta}A=\ct \hat{A}+\st \hat{B}, \quad
\II_{\theta}B=\st \hat{A}-\ct \hat{B}, \quad
\II_{\theta}C=\hat{C}, \quad \II_{\theta}D=-\hat{D}, \\
\II_{\theta}\hat{A}=-\ct A-\st B, \quad
\II_{\theta}\hat{B}=-\st A+\ct B, \quad \II_{\theta}\hat{C}=-C,
\quad \II_{\theta}\hat{D}=D,
\end{gather*}
and
\[\JJ A=B,\quad \JJ C=D,\quad \JJ \hat{A}=\hat{B},\quad \JJ \hat{C}=\hat{D}.\]

Considering $\{J,E_{\theta}'\}$, we have the hypercomplex structure
$\{\II_{\theta}',\,\JJ\}$, with $\JJ$ as above and
\begin{gather*}
\II_{\theta}'A=\ct \hat{A}+\st \hat{B}-\st
\hat{C}+(1+\ct)\hat{D}, \\ \II_{\theta}'B=\st \hat{A}-\ct
\hat{B}+(1+\ct)\hat{C}+\st \hat{D},\\ \II_{\theta}'C=\hat{C},
\quad \II_{\theta}'D=-\hat{D},\\ \II_{\theta}'\hat{A}=-\ct A-\st
B+\st C-(1+\ct)D,\\ \II_{\theta}'\hat{B}=-\st A+\ct
B-(1+\ct)C-\st D,\\ \II_{\theta}'\hat{C}=-C, \quad
\II_{\theta}'\hat{D}=D.
\end{gather*}

For $\{J,E_{\theta}''\}$, we have the hypercomplex structure
$\{\II_{\theta}'',\,\JJ\}$, where $\JJ$ is as above and
\begin{gather*}
\II_{\theta}''A=\ct \hat{A}+\st \hat{B},\quad
\II_{\theta}''B=\st \hat{A}-\ct \hat{B}, \\
\II_{\theta}''C=-\st \hat{A}+(1+\ct)\hat{B}+\hat{C},\quad
\II_{\theta}''D=(1+\ct)\hat{A}+\st \hat{B}-\hat{D}, \\
\II_{\theta}''\hat{A}=-\ct A-\st B,\quad
\II_{\theta}''\hat{B}=-\st A+\ct B, \\
\II_{\theta}''\hat{C}=\st A-(1+\ct)B-C,\quad
\II_{\theta}''\hat{D}=-(1+\ct)A-\st B+D,
\end{gather*}

Finally, for the complex product structure $\{J,\tilde{E}\}$ we have the
hypercomplex structure $\{\tilde{\II},\,\JJ\}$, with $\JJ$ as above and
\begin{gather*}
\tilde{\II}A=-\hat{A}-2\hat{C}, \quad \tilde{\II}B=\hat{B}+2\hat{D}, \quad
\tilde{\II}C=\hat{C}, \quad \tilde{\II}D=-\hat{D},\\ \tilde{\II}\hat{A}=A+2C,
\quad \tilde{\II}\hat{B}=-B-2D, \quad \tilde{\II}\hat{C}=-C, \quad
\tilde{\II}\hat{D}=D.
\end{gather*}
\end{exam}

\vs

\vs
\vs

\n{CIEM, FaMAF, Universidad Nacional de C\'ordoba, Ciudad Universitaria,
(5000) C\'ordoba, Argentina}\\\texttt{andrada@mate.uncor.edu}

\vs 

\n{Dipartimento di Matematica, Politecnico di Torino, Corso Duca degli Abruzzi
24, 10129 Torino, Italy}\\\texttt{salamon@calvino.polito.it}

\enddocument